\documentclass[preprint,12pt]{elsarticle}

\usepackage{setspace,caption}
\usepackage{lipsum}

\usepackage{latexsym,amsmath,amsfonts,amscd}
\usepackage{mathrsfs} 
\usepackage{amsthm}
\usepackage{amssymb}
\usepackage{graphicx}
\usepackage{epstopdf}
\usepackage{ifthen}
\usepackage{etoolbox}
\usepackage{color}
\usepackage{cancel}

\usepackage[usenames,dvipsnames]{xcolor}
\usepackage{soul}
\usepackage[normalem]{ulem}

\usepackage{subcaption} 

\usepackage{stmaryrd} 

\usepackage{nameref}
\usepackage{zref-xr}
\zxrsetup{toltxlabel}
\zexternaldocument*[BilayerPearl-]{../../BilayerPearling/PearlingBilayerV35/PearlingforBilayerV35}

\usepackage{hyperref}

\numberwithin{equation}{section}

\newtheorem{theo}{Theorem}[section]

\newtheorem{definition}[theo]{Definition}

\usepackage{wrapfig}

\usepackage{cancel}

\date{\today}

\def\beal{\begin{align*}}
\def\eeal{\end{align*}}

\def\beq{\begin{equation}}
\def\eeq{\end{equation}}
\def\beqa{\begin{eqnarray}}
\def\eeqa{\end{eqnarray}}
\def\bpm{\begin{pmatrix}}
\def\epm{\end{pmatrix}}

\def\cF{{\cal F}}
\def\cG{{\cal G}}

\def\cQ{{\cal Q}}
\def\cS{{\cal S}}

\def\cZ{{\cal Z}}
\def\eps{\varepsilon}

\def\tu{\tilde{u}}

\def\vV{{\bf V}}

\def\bn{{\bf n}}
\def\bm{{\bf m}}

\def\bN{{\bf N}}
\def\bT{{\bf T}}
\def\bV{{\bf V}}

\def\mbbR{\mathbb{R}}
\def\RR{\mathbb{R}}

\def\NN{\mathbb{N}}
\def\fil{{{\it f}}}
\def\intreach{\int_{-{\sc l}/\eps}^{{\sc l}/\eps}}

\newcommand{\leavethisout}[1] {}

\newcommand{\norm}[1]{\left\lVert#1\right\rVert}
\newcommand{\tck}[1]{{1}}

\newtheorem{asmp}{Assumption}

\newcommand{\draft}[1]{}

\hypersetup{
	colorlinks=true,
	linkcolor=WildStrawberry,
	citecolor=ProcessBlue,
}

\numberwithin{equation}{section}

\newcommand\bt[1]{\begin{tabular}{#1}}
\newcommand\et{\end{tabular}}

\def\bc{\begin{center}}
\def\ec{\end{center}}

\def\vkappa{\vec{\kappa}}
\def\vz{{\bf z}}
\def\Hbl{H}

\def\RR{\mathbb{R}}

\usepackage{lineno}

\journal{Physica D}

\begin{document}

\begin{frontmatter}

\title{Competition and Complexity in Amphiphilic Polymer Morphology}

\author[msu-cm]{Andrew Christlieb}
\ead{andrewch@math.msu.edu}
\address[msu-cm]{Department of Computational Mathematics, Science, and Engineering, Michigan State University}

\author[uofu]{Noa Kraitzman \corref{cor1}}
\ead{noa@math.utah.edu}
\address[uofu]{Department of Mathematics, University of Utah}
\cortext[cor1]{Corresponding author}

\author[msu]{Keith Promislow}
\ead{kpromisl@math.msu.edu}
\address[msu]{Department of Mathematics, Michigan State University}

\begin{abstract}
The strong Functionalized Cahn Hilliard equation models self assembly of amphiphilic polymers in solvent. It supports codimension one and two structures that each
admit two classes of bifurcations: pearling, a short-wavelength in-plane modulation of interfacial width, and meandering, a long-wavelength
instability that induces a transition to curve-lengthening flow. These two potential instabilities afford distinctive routes to changes in codimension and creation of
non-codimensional defects such as end caps and $Y$-junctions.  Prior work has characterized the onset of pearling, showing that it couples strongly to the
spatially constant, temporally dynamic, bulk value of the chemical potential. We present a multiscale analysis of the competitive evolution of codimension one and
two structures of amphiphilic polymers within the $H^{-1}$ gradient flow of the strong Functionalized Cahn Hilliard equation. Specifically we show
that structures of each codimension transition from a curve lengthening to a curve shortening flow as the chemical potential falls through a corresponding critical value.
The differences in these critical values quantify the competition between the morphologies of differing codimension for the amphiphilic polymer mass.
We present a bifurcation diagram for the morphological competition and compare our results quantitatively to simulations of the full system and qualitatively to simulations
of self-consistent mean field models and laboratory experiments. In particular we propose that the experimentally observed onset of morphological complexity arises from a
transient passage through pearling instability while the associated flow is in the curve lengthening regime.
\end{abstract}


\end{frontmatter}

\section{Introduction}\label{Sec:FCHFE}
Amphiphilic molecules are increasingly important in synthetic chemistry where they permit molecular level control of the self assembly of
materials with desirable ionic and electronic conduction properties, \cite{blanazs2009self}. A molecule is amphiphilic with respect to a solvent if
it is comprised of two components, one of which has an energetically favorable or hydrophilic interaction with the solvent and the other with an energetically unfavorable
or hydrophobic interaction. There is a growing literature for the construction and characterization of amphiphilic diblock polymers comprised of polymer chains formed
from adjustable lengths of hydrophilic and hydrophobic polymers covalently bonded together.  Amphiphilic molecules are typically characterized by the Flory-Huggins parameters
indicating the strength of the hydrophobic/hydrophilic interactions and by the aspect ratio of the two components, \cite{fraaije2003model}.

\begin{figure}[h!]
\begin{center}
\begin{tabular}{p{4.0in}p{2.0in}}
\\[0.3in]
    \includegraphics[width=3.7in]{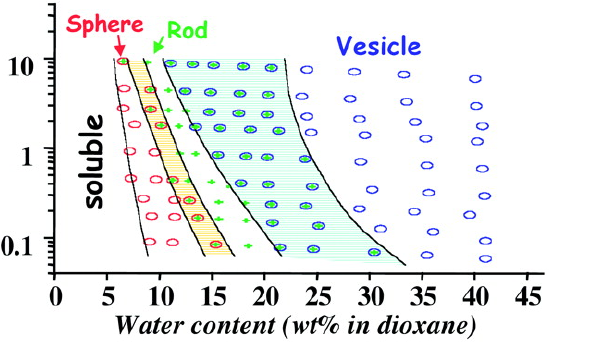} & \\
  \et
  \vspace{-2.7in}
 \begin{tabular}{p{3.7in}c}
 &  \includegraphics[width=2.18in]{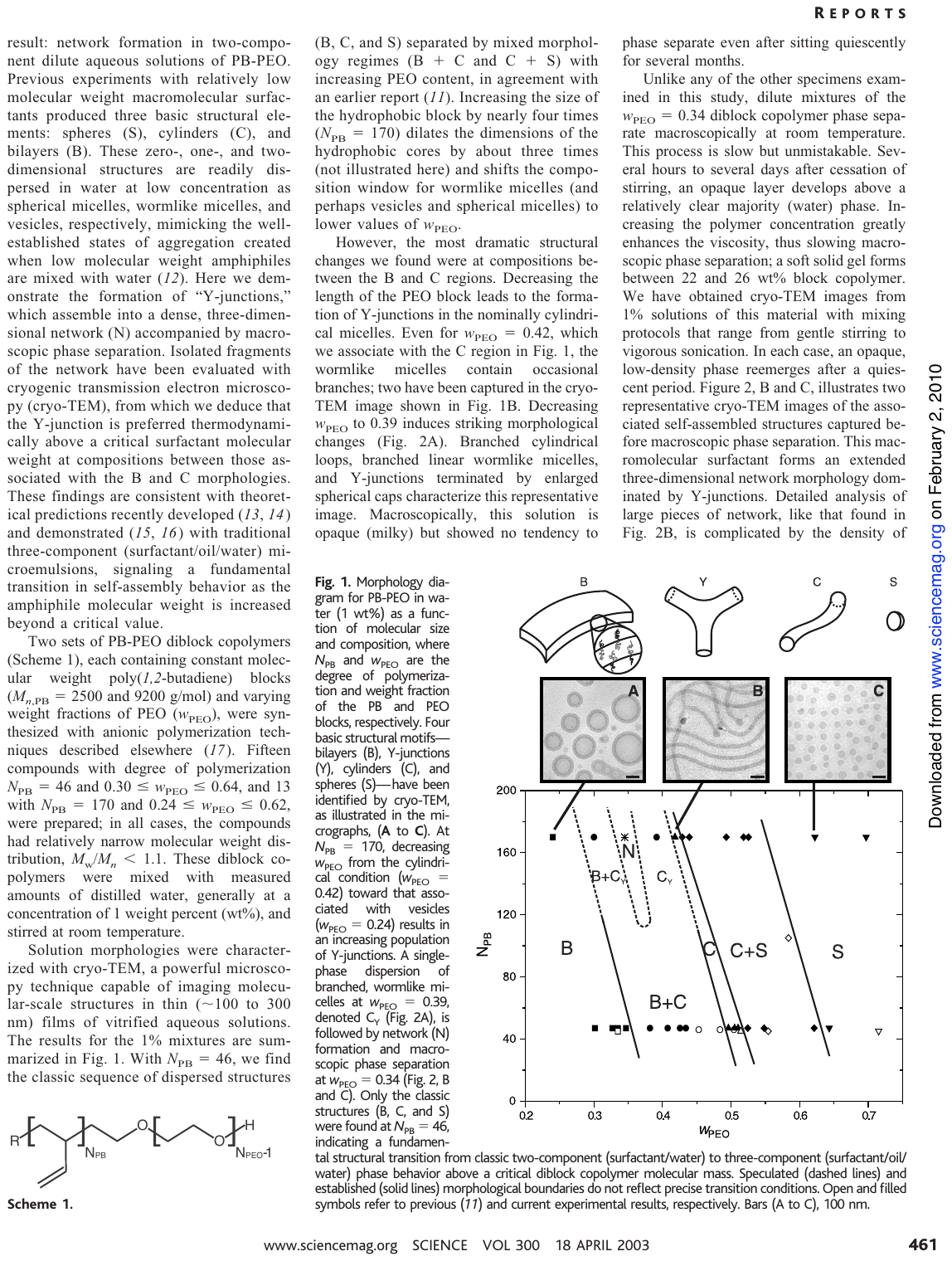}
\end{tabular}
 \vskip -0.1in
  \caption{Morphological bifurcation diagrams for two classes of amphiphilic polymers in solvent. (left) Morphology of Polystyrene (PS) - Polyacrylic acid (PAA) diblocks as function of
 increasing water content in water-dioxane solvent blend (horizontal axis) and polymer \%-weight fraction of the overall mixture (vertical axis).
 Increased water volume fraction drives bifurcation to lower codimensional morphologies without inducing pearling. From \cite{discher2002polymer}. Reprinted with permission from AAAS. (right)
Morphology of amphiphilic Polyethylene oxide (PEO) - Polybutadiene (PB) diblock suspension as function of the PEO weight fraction, $w_{\rm PEO}$ (horizontal axis)
for molecular weights of PB fixed at $N_{\rm PB}=45$ and $170$ (vertical axis). Defect loaded phases are observed for the higher value of $N_{\rm PB}$.
Images of polymer morphologies are included inset, as well as schematic representations of codimension one, two, and three morphologies and `Y' junctions. From \cite{jain2003origins}. Reprinted with permission from AAAS.
}
\label{f:Dicher}
  \end{center}
\end{figure}

When immersed in solvent, amphiphilic polymers self assemble into a wide variety of structures with diverse morphologies
that include codimension one bilayers (hollow vesicles), codimension two filaments (solid rods or cylinders), codimension three micelles (solid spheres), and
various  defect structures with no well defined codimension such as end-caps, ``Y'' junctions, and mixed morphologies. 
The bifurcation structure of these mixtures has been experimentally investigated for a variety of different diblock structures.
Two seminal experiments, presented in Figure\,\ref{f:Dicher}, have partially unfolded this bifurcation structure by studying the impact of varying the
solvent blend to modify strength of the amphiphilic interaction, \cite{discher2002polymer}, and of varying the aspect ratio of amphiphilic polymer, \cite{jain2003origins}.
Increasing the strength of the amphiphilic interaction and decreasing the aspect ratio of the minority phase produce similar results: a sequence
of bifurcations in which codimension-one bilayers yield to codimension-two filaments which yield in turn to codimension-three micelles. In some experiments these
codimensional structures are reported to coexist for a range of parameter values while others show regions of  ``morphological complexity''  characterized by
an abundance of defects.

In this paper we analyze the coexistence, bifurcation, and longtime evolution of well-separated, defect free, codimension one and two structures within the
context of the $H^{-1}$ gradient flows of the strong scaling of the Functionalized Cahn Hilliard free energy.  We use multiscale analysis to derive the curvature
driven evolution of codimension one and codimension two structures that are sufficiently far from self intersection.  We show that the
evolution of morphologies of distinct codimension couples through their exchange of amphiphilic molecules with the bulk. The bulk chemical potential is
spatially constant and varies temporally on a slow time scale. Most significantly, codimension one and two structures switch between a  regularized curve lengthening
and a  curve shortening evolution as the bulk chemical potential passes through critical values. This dichotomy presents a mechanism for morphological competition in
which structures of one codimension grow at the expense of the other. This suggests that, in the absence of defects, the coexistence of structures with distinct codimensions
is not generic however the resolution into  structures of homogeneous codimension  may require a substantial transient. This finding is supported by experimental results,
see \cite{jain2004consequences}, which report that transient structures can persist for months.

There are two mechanisms for a codimension one or two structure to develop an initial defect: self intersection and pearling.
Pearling bifurcations are high-frequency tangential modulations of the width of the codimensional structure. In a companion
paper, \cite{NK-KP-18}, we characterized the onset of the pearling bifurcation, showing that within the strong scaling of the functionalized Cahn Hilliard free
energy the onset of pearling is independent of the shape of the codimension one or two structure, but couples to the transient value of the bulk chemical potential.
Self intersections can be non-local, arising when the initial intersection points are well separated in distance measured along the curve, as in a figure-8 intersection.
They can also be local, as arises when a surface develops large curvatures, such as when the radius of a sphere tends to zero. In both cases the self intersections arise from
the evolution of  the underlying manifold that characterizes the structure. Curve-shortening flows render the manifold smaller and drive their curvatures towards
constant values. Conversely, curve lengthening flows act like backward heat equations for the interfacial curvature and are well known to be locally ill-posed.
We derive a regularized curve lengthening flow that includes a higher-order surface diffusion. The regularized curve lengthening is locally well posed evolution, but causes the
underlying curve to grow and ``meander'' or buckle, and typically leads to finite-time non-local self intersections. In both cases the finite-time singularities can be arrested by the
quenching element of the flow which slows the normal velocity as the far-field chemical potential approaches an equilibrium value. Depending upon initial conditions, our results
can predict either a relaxation to an equilibrium state or provide the mechanisms for the generation of defect states.

Our analysis is particularly relevant to the study of morphologies derived from casting processes in which an initial
suspension of amphiphilic molecules, reflecting a high bulk chemical potential, nucleates out structures of various  codimension which
initially grow, absorbing the amphiphilic molecules from the bulk suspension and lowering the bulk chemical potential. As the chemical
potential falls it may trigger or inhibit pearling bifurcations in one or both codimensions, or trigger transitions from the regularized curve lengthening to
curve shortening.
We compare our asymptotic results to simulations of the full system,
to simulations of self-consistent mean-field density functional models of amphiphilic polymer melts, and to experimental bifurcation diagrams.

\subsection{The Functionalized Cahn Hilliard free energy}

The Functionalized Cahn Hilliard (FCH) free energy models the free energy of a binary mixture of an amphiphilic molecule and a solvent.
It supports stable network morphologies including codimension one bilayers and codimension two filaments as well as pearled morphologies and the defects such as end-caps and $Y$-junctions, \cite{promislow2009pem,gavish2011curvature, dai2013geometric, doelman2014meander, promislow2015existence}.
The FCH free energy takes the form
\beq\label{e:FCHFE}
 \cF(u) := \int\limits_\Omega \frac {1}{2} \left(\eps^2\Delta u-W^\prime(u) \right)^2-
 \eps^p \left(\frac{\eps^2\eta_1}{2}|\nabla u|^2+\eta_2W(u)\right)\, dx,
 \eeq
where $W$  is a smooth double-well potential with local minima at $u=b_\pm$ with $b_-<b_+$. The two wells have unequal depths that are normalized so that
$W(b_-)=0>W(b_+)$ and the left well is non-degenerate in the sense that $\alpha_-:=W''(b_-)>0$. The value of $\alpha_-$ is a key parameter that controls
the rate of exchange of amphiphilic molecules between the bulk and the various morphologies. Here $\eps\ll1$ is small parameter corresponding to the ratio of
length of the amphiphilic molecule to the domain size, $u=b_-$ is associated to a bulk solvent phase, while the quantity $u-b_->0$
is proportional to the density of the amphiphilic phase.  The first term in the integrand of (\ref{e:FCHFE}) is called the Willmore or the quadratic term,
as it denotes the square of a variational derivative of a Cahn Hilliard type free energy. The quadratic term is positive, and we refer to the class of
$u\in H^2(\Omega)$ for which the residual of the quadratic term is small compared to $\eps$ as \emph{morphologies.}  The second grouping of terms in the integrand,
multiplied by $\eps^p$, are called the functionalization terms. The strong and weak scalings of the FCH free energy correspond to the
choice $p=1$ and $p=2$, respectively in (\ref{e:FCHFE}) and represent two natural choices of distinguished limits between the residual
of the quadratic term and the typical scaling of the functionalization terms. In the strong scaling of the FCH, the $O(\eps)$ functionalization terms
typically dominate the generically $O(\eps^2)$ residuals from the quadratic terms, in the weak scaling both terms balance at
$O(\eps^2)$. The analysis of this paper focuses on the {\em strong} scaling of the FCH free energy for which the bifurcation analysis is
more accessible.

The functionalization parameters $\eta_1$ and $\eta_2$ characterize key properties of the amphiphilic molecules.
Specifically $\eta_1>0$ models the strength of the hydrophilic interaction, modeling the propensity of amphiphilic molecules to form monolayers by rewarding
increases in interfacial area or curve length with a decrease in free energy. The parameter $\eta_2\in\RR$ encodes the aspect ratio of the amphiphilic molecule,
as discussed in section\,\ref{s:ABD-sim}. Equivalently these parameters are analogous to the surface and volume energies typical of models of charged solutes
in confined domains, see \cite{scherlis2006unified} and particularly equation (67) of \cite{andreussi2012revised}. With these parameter choices the minus sign in front of
the functionalization terms has great significance --  it incorporates the propensity of the amphiphilic surfactant phase to drive the creation of interface. Indeed, experimental
tuning of solvent quality identifies molecular level phase separation and self assembly in amphiphilic mixtures with the onset of negative values of
surface tension in mesoscale agglomerates, \cite{zhu2009tuning} and \cite{zhu2012interfacial}.

Prior work on FCH gradient flows has focused on the weak scaling, corresponding to $p=2$ in (\ref{e:FCHFE}).
In \cite{dai2013geometric} the authors derived the geometric evolution of bilayers at the quenched mean-curvature flow on the $O(\eps^{-1})$ time scale and as a surface area preserving Willmore flow on the
$O(\eps^{-2})$  time scale. The geometric evolution of codimension two structures was derived in \cite{Dai2015competitive}, obtaining a
curvature driven competitive geometric evolution of the filament curve on the $O(\eps^{-1})$ time scale and a length-preserving Willmore type flow on
the $O(\eps^{-2})$ time scale. Moreover, it was found that the codimension one and two structures can co-exist on the faster $O(\eps^{-1})$
time scale in the weak FCH, but compete on the longer $O(\eps^{-2})$ time scale. However, rigorous investigation of the pearling bifurcation in
the weak FCH is complicated by its leading order coupling to the curvature of the underlying curves. Conversely in the strong scaling of
the FCH, the  pearling bifurcation is independent of morphology and was rigorously characterized in \cite{NK-KP-18} for a wide class of  codimension one and two
structures.


In the remainder of this paper we consider the strong scaling of the FCH free energy. Fixing $\Omega=[0,L]^d\subset{\RR}^d$ for $d=2, 3, \ldots$ and applying
periodic boundary conditions to $H^4(\Omega)$, the first variation of $\cF$, also called the chemical potential~$\mu,$ associated to
a spatial distribution $u\in H^4(\Omega)$ takes the form
\begin{align}\label{e:mu}
	\mu:=\frac{\delta \cF}{\delta u}(u) &= (\eps^2\Delta-W''(u)+\eps\eta_1)(\eps^2\Delta u-W'(u))+\eps\eta_dW'(u),
\end{align}
where~$\eta_d:=\eta_1-\eta_2$. The Functionalized Cahn Hilliard equation
is the associated $H^{-1}$ gradient flow,
\beq
\label{e:FCH-eq}
u_t = \Delta \mu(u),
\eeq
supplemented with periodic boundary conditions on $\Omega$. The choice of the $H^{-1}$ gradient is a reflection of its status as the simplest local gradient that preserves mass.
The mathematical focus of the paper is on the multiscale analysis of the evolution codimension one and two structures. On the~$O(\eps^{-1})$ time scale we find that the $H^{-1}$ gradient flow drives
well separated filament and bilayer structures through a competitive, mean-curvature driven flow mediated through the
common value of the spatially constant far-field chemical potential, $\mu_1$  defined in (\ref{MB-eq:Mu1}).
We show in section 2 that the nonlocal Mullins-Sekerka problem familiar to Cahn Hilliard evolution is present but is unforced, and on the long
time scales we consider the far-field chemical potential relaxes to a trivial, spatially-constant for both codimension one bilayer and codimension
two filament morphologies. As a consequence, the geometric evolution is local.

While spatially constant, the far-field chemical potential $\mu_1=\mu_1(t)$, is temporally dynamic and is linearly proportional to the density of free amphiphilic molecules
 in the bulk. It serves as a key bifurcation parameter, triggering two potential types of instability for each codimension of morphology. Indeed, in \cite{hayrapetyan2015spectra} it is shown for the FCH free energy that the
pearling and self intersection via geometric motion are the only possible mechanisms to generate defects in bilayers.
In the companion paper, \cite{NK-KP-18}, a sharp condition for pearling stability is derived that relates the chemical
potential to the parameter $\eta_d$ and constants that depend implicitly on the
form of the double-well $W$. Specifically the bilayers are stable with respect to the pearling bifurcation if and only if
\beq
    \mu_1S_b+\eta_d\lambda_{b,0}\norm{\psi_{b,0}}_2^2<0,
\label{eq:B-PS}
\eeq
and similarly filaments are pearling stable if and only if
\beq
    \mu_1S_\fil+\eta_d\left(\norm{\psi_{\fil,0,0}'}_{L_R}^2+\lambda_{\fil,0,0}\norm{\psi_{\fil,0,0}}_{L_R}^2\right)<0,
\label{eq:F-PS}
\eeq
where~$\lambda_{b,0}$ is the ground-state eigenvalue of the linear operator~$L_{b,0}$ , defined in~(\ref{e:L0-def}), with eigenfunction~$\psi_{b,0}$,  and $\lambda_{\fil,0,0}$ is the ground state eigenvalue of the linear operator~$L_{\fil,0,0}$, defined in~(\ref{CS-eq:Lm}), with the corresponding eigenfunction~$\psi_{\fil,0,0}$.
The constants $S_b,S_\fil$ are the {\bf shape factors} of the bilayers and the filaments, respectively, defined by the relations
\begin{align} \label{eq:SF-def}
    S_b:=\int_\RR\Phi_{b,1} W'''(\phi_b)\psi_{b,0}^2\,dz,\quad
    S_\fil:=2\pi\int_0^\infty \Phi_{\fil,1}W'''(\phi_\fil)\psi_{\fil,0,0}^2\,RdR.
\end{align}
Here $\phi_b$ and $\phi_\fil$ are the bilayer and filament profiles, defined in (\ref{e:BLcp}) and (\ref{e:FLcp}), while
$\Phi_{b,1}$ and $\Phi_{\fil,1}$, defined in (\ref{e:Phi_bj-def}) and (\ref{e:Phi_fj-def}), encode the impact of
a change in chemical potential on the shape of the bilayer and filament, respectively.  For each codimension, if the shape factor is negative then pearling stability is favored
by large (positive) values of $\mu_1$, while if it is positive then pearling stability is favored for small (negative) values of $\mu_1$. In \cite{promislow2015existence} the existence of
pearled codimension one circular and flat equilibrium was demonstrated in $\RR^2$ for the strong FCH.

\begin{figure}[h!]
\begin{center}
\begin{tabular}{c}
    \includegraphics[width=6.0in]{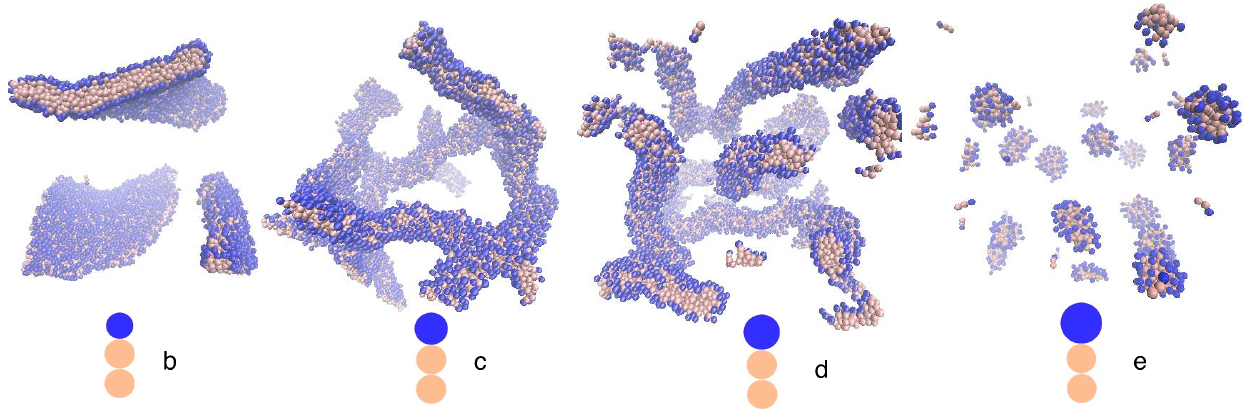}
\end{tabular}
 \vskip -0.1in
  \caption{Molecular dynamics simulations of amphiphilic molecules with distinct aspect ratios. Increasing the size of the hydrophilic head group, relative to a fixed tail leads to a preference for morphologies with
  increasing codimensions: (b) bilayers with edge caps (codim 1), (c)  branched filaments (codim 2), (d) unbranched filaments (codim 2), and (e) micelles (codim 3). Morphologies with higher codimension have a higher density of free-floating amphiphilic molecules in the far-field (bulk), corresponding to a higher critical value of $\mu_1$. Reprinted from \cite{COOKE2006COUPLING}, with permission from Elsevier.
}
\label{f:Aspect}
  \end{center}
\end{figure}

\subsection{Summary of Analytical Results}
Our main analytical result is the derivation of curvature driven flow laws valid on slow time $\tau=\eps t$ for codimension one and two
structures embedded in $\Omega\subset\RR^3$ via a multiscale analysis. Bilayer and filament morphologies are defined as dressings of collections of
admissible interfaces and curves, respectively. We consider an admissible codimension one interface $\Gamma_b$, see Definition\,\ref{def:admissible}, with $O(1)$
area and codimension two curve $\Gamma_\fil$, see Definition\,\ref{def:admissible_f}, with $O(\eps^{-1})$ length. We assume the reaches of these two classes of surfaces, as defined in
(\ref{e:whiskNbd}) and (\ref{e:whiskNbd-f}) respectively, are all disjoint and introduce the composite solution defined in (\ref{e:Composite}). The evolution
of these composite solutions under the gradient flow (\ref{e:FCH-eq}) is parameterized at leading order by the triple $(\Gamma_b, \Gamma_\fil, \mu_1)$  according to
their $\eps$-scaled normal velocities
\begin{align}
	V_b & = \nu_b(\mu_1-\mu_b^*)\Hbl_0+\eps k_b \Delta_s \Hbl_0, \label{e:NVb}\\
	\bV_\fil &= -\left[\nu_\fil(\mu_1-\mu_\fil^*)\vec{\kappa}+\eps k_f \partial_s^2 \vec{\kappa}\right]. \label{e:NVf}
\end{align}
Here we have introduced the mean curvature $H_0$ and Laplace-Beltrami operator $\Delta_s$ of $\Gamma_b$ and vector curvature $\vec{\kappa}=(\kappa_1,\kappa_2)$
and surface diffusion $\partial_s^2$ of $\Gamma_\fil$. On this slow time $\tau=\eps t$, the chemical potential satisfies
\begin{align}\label{e:Mu1Evol}
    \frac{d\mu_1}{d\tau} = -\frac{\alpha_-^2}{|\Omega|}\Bigg[&m_b\int_{\Gamma_b}\nu_b\left({\mu_1}-\mu_b^*\right)\Hbl_0^2-\eps k_b |\nabla_s H_0|^2\,ds \\
    &+
     2\pi m_\fil \eps \int_{\Gamma_\fil}
    \nu_\fil\left({\mu_1}-\mu_\fil^*\right)|\vec{\kappa}|^2-\eps|\partial_s\vec{\kappa}_\fil|^2\,ds\Bigg]
    +O(\eps^2)\nonumber,
\end{align}
where the constants $\nu_b>0$, $k_b>0$, and $m_b>0$, are defined in (\ref{e:MB-const}) 
while $\nu_\fil>0$, $k_\fil>0$, and $m_\fil>0$  defined in (\ref{MP-eq:S1}). While the surface diffusion terms are formally lower order, they
are leading order singular perturbations that keep the resulting flow locally well-posed.

 The system (\ref{e:NVb})-(\ref{e:Mu1Evol}) describes the competitive dynamics between  codimension one and two morphologies. The key critical values
\begin{align}
\mu_b^* &= -\frac{k_b}{2\nu_b}(\eta_1+\eta_2), \label{e:nub-def} \\
\mu_\fil^* & = \frac{k_\fil}{\nu_f} \eta_1, \label{e:nuf-def}
\end{align}
indicate the value of $\mu_1$ at which the rates of amphiphilic molecule insertion and ejection are balanced for bilayers and filaments, respectively.
Specifically, if the far-field chemical potential lies above this number, then structures of the corresponding codimension  will grow. Indeed, the rate of change of
the area of bilayers is given by
\beq
\label{e:Gb-grow}
\partial_\tau |\Gamma_b| = \int_{\Gamma_b} V_b H_0\,ds = \int_{\Gamma_b} \nu_b (\mu_1-\mu_b^*) H_0^2 -\eps k_b |\nabla_s H_0|^2\,ds,
\eeq
with the corresponding expression
\beq
\label{e:Gf-grow}
\partial_\tau |\Gamma_\fil| = \int_{\Gamma_\fil} \bV_\fil\cdot\vec{\kappa}\, ds = \int_{\Gamma_\fil} \nu_b (\mu_1-\mu_\fil^*) |\vec{\kappa}|^2 -\eps k_\fil |\partial_s \vec{\kappa}|^2\,ds,
\eeq
for the length of the filaments. The competitive dynamics system provides the leading order evolution so long as the interfaces remain admissible with disjoint reaches and the
$\mu_1$-dependent pearling conditions (\ref{eq:B-PS})-(\ref{eq:F-PS}) hold. The surface diffusion terms in (\ref{e:NVb})-(\ref{e:NVf}) are relevant to mass balance only if the gradients of
the curvatures become asymptotically large or if $\mu_1$ becomes asymptotically close to one of the critical values $\mu_b^*$ or $\mu_\fil^*$. In particular, it follows from (\ref{e:Mu1Evol}) that
net growth of bilayers and filaments corresponds to a decrease in $\mu_1$, while large curvature gradients enhance the ejection rates and increase the value of $\mu_1$.

{\bf Remark}:
For the codimension two filament term to contribute to the evolution of the chemical potential, $\mu_1$, at leading order, we assume that their collective length
 $|\Gamma_\fil|$ is $O(\eps^{-1})$. Our generic assumption on the codimension one phase is the surface area $|\Gamma_b|=O(1)$, so that
both bilayers and filaments occupy an $O(\eps)$ volume fraction. This limits the applicability of the asymptotic results as the assumption of disjoint reaches becomes
non-generic, and represents a significant caveat in the application of our analytical results. Since the geometric flow reduction does not apply to morphologies with defects,
both a large number of short filaments or a small number of long filaments complicate the non-self intersection assumption.
Our analysis applies to each disjoint component of filament morphology separately, and while $\eps\ll1$ is small, it is viewed as fixed within the model and need not be vanishingly small.
A more detailed analysis of mass scaling in the $\eps\to0$ convergence issues for FCH models can be found in \cite{dai2019constrained}.


For $\mu_1> \mu_b^*$ we call the normal velocity (\ref{e:NVb}) a regularized (codimension one) curve lengthening flow and a (codimension one) curve shortening flow if $\mu_1<\mu_b^*$,
with similar terminology for the codimension two flow based upon the sign of $\mu_1-\mu_\fil^*.$
When the structures have a homogeneous codimension, then in the absence of singularities in the curvature flow, equation (\ref{e:Mu1Evol}) drives the chemical potential $\mu_1$
to the corresponding critical value, $\mu_b^*$ or $\mu_\fil^*$, and the leading-order term in the geometric flow goes to zero, and the system is said to be ``quenched''.
To illustrate the nature of the geometric flow, it is instructive to rewrite it as a corresponding evolution equation for the curvatures. For codimension one structures
in two space dimensions, up to tangential reparameterization it takes the simple form
\beq
 \partial_\tau H_0 = -(\Delta_s +\Hbl_0^2) V_b = -(\Delta_s +\Hbl_0^2)\left(\nu_b(\mu_1-\mu_b^*)\Hbl_0+\eps k_b \Delta_s \Hbl_0\right).
 \eeq
 For the curve lengthening regime, the dominant term is a backward heat equation, with a fourth-order regularization and a $\Hbl_0^3$ nonlinearity
 with a negative (stable) coefficient. For the curve shortening regime, both differential terms are stabilizing, but the cubic nonlinearity has a positive coefficient
 that supports finite time singularity which may be arrested by the fourth order regularization.

In numerical simulations the curve lengthening flows, with $\mu_1>\mu^*_{b,\fil}$ show distinct regimes. For positive but $O(\eps)$ values
of $\mu_1-\mu^*_b$ the bilayer interfaces will bend and buckle at $O(1)$ length scales, leading to shapes reminiscent of a meandering river.
This regime is called a ``meandering flow'', and is studied rigorously in \cite{CP-19}. For $O(1)$ positive values of $\mu_1-\mu_b^*$ the curve lengthening flows
can lead to growth of high-curvature regions which self intersect on a $\tau=O(1)$ length scale. For filaments this can lead to the formation of many closed loops,
see Figures\,2 and 3 of  \cite{jain2003origins} for experimental examples or Figure\,\ref{f:Transman} of this article for an example of meandering motion within the FCH gradient flow.

\subsection{Competition and Morphological Complexity}

Our main scientific results are the conjectured mechanisms for the morphological bifurcations observed in the casting of amphiphilic suspensions.
In particular we propose a mechanism for the onset of the so-called ``morphological complexity'' observed in the experimental casting processes
presented in Figure\,\ref{f:Dicher} (right). In a casting process amphiphilic molecules are dispersed (mixed) in a solvent and the mixture is allowed to relax, generically leading to self-assembly of structures with
distinct codimension.  For the shorter chains,  $N_{PB}=45$ corresponding to the lower horizontal row of symbols, casts from molecules with increasing weight fraction
of the amphiphilic PEO component ($w_{\rm PEO}$ -- horizontal axis)  lead to a series of codimensional bifurcations in which the self-assembly prefers structures with
increasing codimension.  Indeed experiments produce first bilayers (codimension one) marked (B), then bilayers coexisting with cylinders (codimension two)
marked (B+C), then cylinders, then cylinders coexisting with spheres (codimension three) marked (S), and finally spheres. A similar series of codimensional bifurcations are presented in
Figure\,\ref{f:Dicher}(left) and Figure\,\ref{f:Aspect}. In the former the bifurcations are induced by lowering the percentage of water in the water-dioxane solvent which forms
the basis for the casting. In the latter they are realized within coarse-grained molecular dynamics simulations of casting processes by increasing the size of the hydrophilic head group,
and hence the aspect ratio, of a simple three-group amphiphilic molecule.

We conjecture that the competitive dynamics implicit in the system (\ref{e:NVb})-(\ref{e:Mu1Evol}) forms the basis for these codimensional bifurcations.
In a casting process there is initially a relatively large density of dispersed amphiphilic molecules, corresponding to a high value of the far-field chemical potential $\mu_1$.
As various structures self-organize, amphiphilic molecules are removed from the far-field and the scalar value of $\mu_1$ falls. The critical values $\mu_b^*$ and  $\mu_\fil^*$,
given in (\ref{e:nub-def})-(\ref{e:nuf-def}), guage the relative ability of the corresponding bilayer and filament morphologies to absorb and retain amphiphilic molecules
from the far-field (bulk) environment. The morphology with the lowest corresponding critical value will, in the absence of defects, lower the value of $\mu_1$ and drain the
mass of the morphologies with higher critical values of $\mu_1$. Increasing values of either $\eta_1$ or $\eta_2$ will drive $\mu_b$ greater than $\mu_\fil$ and trigger a competitive
imbalance that favors codimension two filaments over codimension one bilayers. We argue in Section 4 that increasing the aspect ratio of the amphiphilic molecule corresponds to an
increase in the value of $\eta_2$, while decreasing the percentage of water within the water-dioxane solvent blend corresponds to an increase in the value of  $\eta_1$. These produce
shifts in $\mu_{b,\fil}^*$ in agreement with bifurcation from codimension one to codimension two. The coexistence of codimension one and two structures for large parameter ranges are
not supported by the analysis. However the time scale to reach equilibrium can be quite long, \cite{jain2004consequences} suggest times on the order of months, and we propose that
longer experimental trials may decrease the size of the regions of coexistence.  We do not present an analysis of codimension three micelles within this work as they do not have a
spatially extended direction that can accommodate incremental growth, rather simulations suggest that they swell, form dumbbell shapes, and then break into distinct micelles. This
behavior is outside the scope of our analysis.

For the longer $N_{PB}=170$ chains in Figure\,\ref{f:Dicher}(right), increasing $w_{\rm PEO}$ one finds that the codimensional bifurcation structure is interrupted by the onset of
so-called ``morphological complexity''. Specifically, the casting sequences yield bilayers, bilayers coexisting with branched filaments,
strongly connected network morphologies,  and $Y$-junction dominated filaments, before reverting to the familiar codimension two and codimension three structures.
The term morphological complexity refers not only to the wide variety of possible outcomes, but also to the difficulty in controlling the outcomes, see
\cite{jain2003origins} and \cite{jain2004consequences}. Our second conjecture is that morphological complexity arises from the interplay between the pearling bifurcation,
the competitive dynamics, and the evolving value of $\mu_1$.  Indeed, the criteria for pearling stability depends upon the value of $\mu_1$, see (\ref{eq:B-PS}) and (\ref{eq:F-PS}),
and as $\mu_1$ deceases during the casting it may trigger or inhibit pearling stability. In particular, in section 4.4.2 we present regimes in which bilayers have a
competitive advantage over filaments, but are {\sl transiently} pearling unstable, while filaments are globally pearling stable. This cascade of bifurcations provides mechanisms
to produce complex blends of defects and morphologies and affords a clear mechanism for hysteresis.  In such an environment the ultimate outcome of a given casting
process could depend sensitively on secondary effects such as the rate at which amphiphilic molecules are initially added to the dispersion or upon spatial inhomogeneities.
The spatial complexity of the end states in this regime is born out both by experiments and by simulations of the FCH free energy, see Figure\,\ref{f:JB}(center and right).

We emphasize that the morphological complexity conjecture encompasses structures with codimensional defects that are outside of our analysis. Moreover,
the pearling bifurcation cannot be robustly suppressed within the FCH gradient flows. The experiments and simulations exhibiting the simple codimensional bifurcation route
do not display signs of pearling bifurcation. In particular there is no mechanism within the FCH energy to explain why pearling would be expressed in longer polymers and
inhibited in shorter but otherwise identical polymers. These limitations of the model are expanded upon in the discussion of Section\,\ref{s:Discussion}. In \cite{promislow2017existence}
and \cite{DPV-19} two-component extensions of the FCH are proposed which posses more detailed internal layer structure and afford precise mechanisms to robustly inhibit
pearling bifurcations.

In Sections 2 and 3, respectively, we apply a multiscale analysis to derive the long time-scale evolution of admissible codimension one and two structures under the FCH equation. In
particular we extract the coupling of the bulk chemical potential on the curvature driven flow. In section 4 we extend these results to composite morphologies,
deriving the laws governing their competitive evolution. We present bifurcation diagrams that show the regions of pearling stability and curve shortening of each morphology, and compare
them to simulations of the FCH equation, to self consistent mean field density functional theory simulations, and to the experimental bifurcation diagrams presented in Figure\,
\ref{f:Dicher}.  These are the first results that explicitly quantify the complexity of the transients associated to amphiphilic polymer blends and identify the role of the pearling bifurcation in the generation of complex network morphologies.

\section{Geometric Evolution of Codimension One Structures}\label{Sec:Bilayer}

In this section we derive the geometric evolution of admissible codimension one interfaces, which we refer to as bilayers. These
calculations are carried out in~$\RR^d$ for $d\geq 2$ but we will restrict our attention to $\RR^3$ for the analysis of filaments.
We consider the local, mass-preserving $H^{-1}$ gradient flow of
the {strong FCH} given in (\ref{e:FCH-eq}). The multiscale analysis of this section follows closely from the calculations of \cite{Dai2015competitive}, which considered the weak
scaling of the FCH gradient flow. For brevity we present only the main calculations.

It is well known that for the single-layer interfaces supported in Cahn-Hilliard type models, that is codimension one interfaces which separate distinct phases, the $O(1)$ and $O(\eps^{-1})$ time
scales yield Stefan and Mullins-Sekerka problems for the interfacial motion \cite{pego1989front, alikakos1999periodic, alikakos2000mullins}.
For single layers motion of the interface requires transport of materials on either side.
For bilayers we derive theses reduced flows, but they have trivial solutions, as the interfacial motion of an interface
with the same material on either side does not require long-range transport but is facilitated by permeation. The result is a local geometric flow, driven by
membrane curvatures and coupled to the bulk value of the chemical potential.

\subsection{Admissible codimension one manifold and their dressings}

Given a smooth, closed $(d-1)$--dimensional manifold $\Gamma_b$ immersed in $\Omega\subset\RR^d$,
we define the local  ``whiskered''  coordinates system in a neighborhood of $\Gamma_b$ via the mapping
\beq
\label{e:codim1-cov}
 x=\rho(s,z):= \zeta_b(s)+\eps \nu(s)z,
\eeq
where $\zeta_b:\cS\mapsto\RR^d$ is a local parameterization of $\Gamma_b$ and $\nu(s)$ is the outward unit normal to $\Gamma_b.$
The variable $z$ is often called the $\eps$-scaled, signed distance to $\Gamma_b$, while the variables $s=(s_1,...,s_{d-1})$ parameterize the tangential
directions of $\Gamma_b$.
\begin{definition}\label{def:admissible}
For any $K, \ell > 0$ the family, $\cG_{K,\ell}^b$,  of admissible interfaces is comprised of closed (compact and without boundary),
oriented $d-1$ dimensional manifolds $\Gamma_b$  embedded in $\mbbR^d$, which are far from self-intersection and with a smooth second fundamental form.  More precisely,
\vskip 0.01in
\begin{tabular}{lp{5.0in}}
 (i)& The  $W^{4,\infty}(\cS)$ norm of the 2nd Fundamental form of $\Gamma_b$ and its principal curvatures are bounded by $K$.\\
 (ii)& The whiskers of length $3\ell < 1/K$, in the unscaled distance, defined for each~$s_0\in \cS$ by, $w_{s_0}:=\{x: s(x)=s_0, |z(x)|<3\ell/\eps\}$,  neither intersect each-other
nor $\partial\Omega$ (except when considering periodic boundary conditions). \\
(iii)& The surface area, $|\Gamma_b|$, of $\Gamma_b$ is bounded by $K$.
\end{tabular}
\end{definition}

For an admissible codimension one interface $\Gamma_b$ the change of variables $x \rightarrow \rho(s,z)$ given by (\ref{e:codim1-cov}) is a $C^4$ diffeomorphism
on the reach of $\Gamma_b$, defined as the set
\beq
\label{e:whiskNbd}
\Gamma_{b,\ell} := \left\{ \rho(s,z) \in \mbbR^d \Bigl| s \in \cS, -\ell/\eps \le z \le \ell/\eps \right\}\subset\Omega,
\eeq
with complement $\tilde{\Gamma}_{b,\ell}:=\Omega\backslash \Gamma_{b,\ell}.$
On the reach we may expand the Cartesian Laplacian in terms of the Laplace-Beltrami operator $\Delta_s$ and the curvatures,
\beq\label{e:LaplacianBl}
 \eps^2\Delta_x 
 =\partial_z^2+\eps H_0(s)\partial_z +\eps^2(zH_1\partial_z +\Delta_s)+O(\eps^3),
\eeq
where $H_i(s)$ is related to the $\text{i}^{\text{th}}$ powers of the curvatures
\beq
    H_i = (-1)^i\sum_{j=1}^{d-1}k_j^{i+1},
\eeq
and, in particular, $H_0$ is the mean curvature of $\Gamma_b$. See \cite{dai2013geometric} for more details.

\begin{definition} Given an admissible codimension one interface $\Gamma_b\in\cG_{K,\ell}^b$ and $f:\RR\rightarrow\RR$ which tends to constant value~$f_{\infty}$ at an
exponential rate as~$z\rightarrow\pm\infty$,  then we define the $H^2(\Omega)$ function
\beq
\label{def:bl-dress}
    f_{\Gamma_b}(x):=f(z(x))\chi(|z(x)|/\ell)+f_\infty(1-\chi(|z(x)|/\ell)),
\eeq
where~$\chi:\RR\rightarrow\RR$ is a fixed, smooth cut-off function which takes values one on~$[0, 1]$ and $0$ on $[2,\infty)$.
We call $ f_{\Gamma_b}\in L^2(\Omega)$ the {\bf dressing} of~$\Gamma_b$ with~$f\in L^2(\RR)$,  and by abuse of notation will drop the $\Gamma_b$ subscript
when doing so creates no confusion.
\end{definition}

Within the reach $\Gamma_{b,\ell}$ of an admissible $\Gamma_b$  the quadratic term within the FCH, (\ref{e:FCHFE}), can be
re-written in the codimension one whiskered coordinates system (\ref{e:codim1-cov}). Setting the quadratic term equal to zero,
and formally taking the leading order terms in $\eps$ leads to a second-order ODE in $z$. The bilayer profile $\phi_b$, is defined to be the
solution of this equation
\beq
\label{e:BLcp}
 \partial_z^2 \phi_b = W'(\phi_b),
 \eeq
which is homoclinic to the left well $b_-$ of $W$. We denote by $U_b\in L^2(\Omega)$ the dressing of $\Gamma_b$ by $\phi_b\in L^2(\RR)$, and introduce the
associated linear operator
 \beq
 \label{e:L0-def}
 L_{b,0}:=\partial_z^2 - W''(\phi_b).
 \eeq
This is a Sturm-Liouville operator on $L^2(\RR)$ and has a positive ground-state eigenvalue $\lambda_{b,0} > 0$ with eigenfunction $\psi_{b,0} \geq 0$ and
a translational eigenvalue $\lambda_{b,1} = 0$ associated to the eigenfunction $\psi_{b,1}=\phi_b'$.
In addition, we define the functions
\beq
\Phi_{b,j}:=L_{b,0}^{-j} 1, \label{e:Phi_bj-def}
\eeq
for $j=1, 2$ which converge to a non-zero value at $z=\pm\infty$. We also define their $\Gamma_{b}$ dressings, which are denoted by $\Phi_{b,1}$ and $\Phi_{b,2}$
by the abuse of notation mentioned above.

\subsection{Inner and outer expansions}
Assuming initial data arising from the dressing of an admissible initial codimension one interface~$\Gamma_b(t_0)\in\cG_{K,\ell}^b$,
we describe the coupled geometric evolution of the interface and the far-field chemical potential as a flow in time~$t$.
We consider formal, multi-scale analysis of the density~$u$ and the chemical potential~$\mu$.
In the far-field or bulk region, ~$\tilde\Gamma_{b,\ell}:=\Omega\backslash \Gamma_{b,\ell}$, the outer solution~$u$ has the expansion
\begin{align}
    \label{MB-eq:uOuterSol}u(x,t) &= u_0(x,t)+\eps u_1(x,t)+O(\eps^2).
\end{align}
Within the reach or inner region, $\Gamma_{b,\ell}$,  we use the whiskered coordinates with the $\eps$-scaled distance $z$ to
the interface. The standard assumption is that the inner solution $\tu$ is smooth in the tangental $s$-variables. However
the leading order result is a curvature driven flow whose coefficient may switch sign. Flow against curvature is not locally well posed
due to uncontrollable growth in  high-frequency growth terms. To regularize this, we incorporate a two-scale tangential expansion,
introducing the fast tangential variable  $S:=s/\sqrt{\eps}$, so that the inner variable admits the expansion.
\begin{align}
    \label{MB-eq:uInnerSol}u(x,t)&=\tilde u(s,S,z,\tau) = \tilde u_0(s,S,z,\tau)+\eps\tilde u_1(s,S,z,\tau)+O(\eps^2).
\end{align}
The inclusion of the fast tangential variable promotes a formally lower order surface diffusion term
to the leading order, where it regularizes the curve lengthening  flow. The normal velocity~$V=V(s,S,t)$ of~$\Gamma_b$  is denoted by
 \beq\label{MB-eq:VDef}
    V(s,S,t):= 
   - \eps\frac{\partial z}{\partial t}.
 \eeq
On the slow time $\tau=\eps t$,  the time derivative of the inner density function~$\tilde u$, defined in~(\ref{MB-eq:uInnerSol}),
combined with the normal velocity,~(\ref{MB-eq:VDef}), takes the form
\beq\label{MB-eq:ut}
    \frac{\partial \tilde u}{\partial t} = -\eps^{-1}V(s)\frac{\partial \tilde u}{\partial z}+\frac{\partial \tilde u}{\partial \tau}\frac{\partial \tau}{\partial t}.
\eeq
At the interface we have the standard matching conditions
\beq\label{MB-eq:MatchCond}
   \lim_{h\rightarrow 0^\pm} u(x+ h\nu,\tau) = \lim_{z\rightarrow \pm\infty}\tilde u(s,S,z,\tau).
\eeq
which reduce to the relations
\begin{align}
    u_0^{\pm}(x,\tau) &= \lim_{z\rightarrow\pm\infty}\tilde u_0(s,S,z,\tau),\label{MB-eq:MatchU0}\\
    u_1^{\pm}(x,\tau)+z\partial_\nu u_0^{\pm}(x,\tau) &= \lim_{z\rightarrow\pm\infty}\tilde u_1(s,S,z,\tau).
\end{align}
where~$\partial_\nu$ is the derivative in the normal direction of~$\Gamma_b$, and~$u_i^{\pm}$ denote the values of the limits of the left-hand side
of (\ref{MB-eq:MatchCond}) as $h\to0^\pm$ respectively.

The chemical potential, $\mu$, defined in (\ref{e:mu}), admits similar outer and inner expansions. The terms  of its outer expansion are slaved
to the density $u$ through the outer relations,
\begin{align}
	\label{MB-eq:Mu0}\mu_0 =& W''(u_0)W'(u_0),\\
	\label{MB-eq:Mu1}\mu_1 =& (W'''(u_0)u_1-\eta_1)W'(u_0)+(W''(u_0))^2u_1+\eta_dW'(u_0).
\end{align}
For the inner expansion we introduce the nonlinear operators $P$ and $Q$ to rewrite the chemical potential (\ref{e:mu}) as
\beq
\tilde\mu = P(\tu) Q(\tu) +\eps \eta_d W'(\tu).
\eeq
In the multiscale tangential variables the Laplacian expansion
(\ref{e:LaplacianBl}) takes the form
\beq\label{e:Lap-ms}
 \eps^2\Delta_x 
 =\partial_z^2+\eps\left( H_0(s,S)\partial_z +\Delta_S\right) +\eps^2(zH_1(s,S)\partial_z + z D_{2,S}+ \Delta_s)+O(\eps^3),
\eeq
where $\Delta_S$ is the scaled Laplace-Beltrami operator and $D_{2,S}$ denotes a higher order elliptic term in $S$. Details on the $D_{2,S}$ term can
be found in section 6 of \cite{hayrapetyan2015spectra} however its precise form is immaterial to our presentation. We combine the expansion of the Laplacian
and the inner solution to obtain an expansion for $ P=P_0+\eps  P_1+\eps^2 P_2+ \ldots$, where
\begin{align}
\label{e:P0}
P_0 = & \partial_z^2 - W''(\tu_0), \\
\label{e:P1}
 P_1 = & H_0(s,S) \partial_z +\Delta_S - W'''(\tu_0)\tu_1 + \eta_1, \\
\label{e:P2}
 P_2 = & z H_1(s,S) \partial_z + D_{2,S} + \Delta_s - W'''(\tu_0)\tu_2 +\frac12 W^{(4)}(\tu_0)\tu_1^2,
\end{align}
and for $Q=Q_0+\eps Q_1+\eps^2 Q_2 +\ldots$
\begin{align}
\label{e:Q0} Q_0 = & \partial_z^2\tu_0 - W'(\tu_0), \\
\label{e:Q1}  Q_1 = & H_0 \partial_z \tu_0 +\Delta_S\tu_0 +(\partial_z^2- W''(\tu_0))\tu_1, \\
\label{e:Q2} Q_2 = & z H_1 \partial_z\tu_0 + D_{2,S}\tu_0 + \Delta_s\tu_0 +(H_0\partial_z+\Delta_S)\tu_1+ (\partial_z^2-W''(\tu_0))\tu_2 -\frac12 W'''(\tu_0)\tu_1^2.
\end{align}
With these reductions  we expand the inner chemical potential as
\begin{align}
\label{MB-eq:tildeMu0}\tilde\mu_0=& P_0Q_0\\
\label{MB-eq:tildeMu1}\tilde\mu_1 =& P_1Q_0+P_0Q_1 +\eta_d W'(\tu_0), \\
\label{MB-eq:tildeMu2}\tilde\mu_2 =& P_0Q_2+P_1Q_1 +P_2Q_0+\eta_d W''(\tu_0)\tu_1 .
\end{align}
The second order form of the $H^{-1}$ gradient induces inner-outer matching conditions for the chemical potential,
\begin{align}
    \mu_0^{\pm}(x,t) &= \lim_{z\rightarrow\pm\infty}\tilde \mu_0(z,s,S,t),\\
    \mu_1^{\pm}(x,t)+z\partial_\nu \mu_0^{\pm}(x,t) &= \lim_{z\rightarrow\pm\infty}\tilde \mu_1(s,S,z,t) \label{MB-eq:MatchMu1}\\
    \mu_2^{\pm}(x,t)+z\partial_\nu \mu_1^{\pm}(x,t)+\frac{1}{2}z^2\partial_\nu^2\mu_0^\pm(x,t)&= \lim_{z\rightarrow\pm\infty}\tilde \mu_2(s,S,z,t).\label{MB-eq:MatchMu2}
\end{align}

\subsection{Time scale~$\tau=\eps t$: quenched curvature driven flow}
We focus on the first relevant slow time-scale, $\tau = \eps t$, inserting the time derivative and chemical potential expansions
into the FCH gradient flow, (\ref{e:FCH-eq}). As the interface $\Gamma_b$ is codimension one, it separates the region $\Omega$ into
two disjoint sets, $\Omega\backslash\Gamma_b=\Omega_+\cup\Omega_-$ with the normal to $\Gamma_b$ pointing towards $\Omega_+$.
In this bulk region we obtain the relations
\begin{align}
O(1): & \hspace{0.1in} 0 =  \Delta \left(W''(u_0)W'(u_0)\right), &\text{in }\Omega_-\cup\Omega_+, \label{MB-eq:t1Outer0}\\
O(\eps): &\hspace{0.1in}  u_{0,\tau} = \Delta \left((W'''(u_0)u_1-\eta_1)W'(u_0)+(W''(u_0))^2u_1+\eta_dW'(u_0)\right),&\text{in }\Omega_-\cup\Omega_+.
\label{MB-eq:t1Outer01}
\end{align}
The relevant solution to the $O(1)$ relation is the spatially constant density, $u_0 = b_-$ for which $\mu_0=0.$ With this reduction and the fact that $W'(b_-)=0$ and
$W''(b_-):=\alpha_->0$, the $O(\eps)$ relation reduces to
\beq\label{MB-eq:t1Outer1}
    0 = \Delta u_1 \quad\text{in }\Omega_-\cup\Omega_+,
\eeq
which is subject to interior layer matching and exterior boundary conditions derived in the sequel.

In the inner region we supplement the expansion (\ref{MB-eq:uInnerSol}) with the form of the Laplacian in inner variables,
given in (\ref{e:LaplacianBl}). Collecting orders of $\eps$ we find
\begin{align}
	O(\eps^{-2}):&\quad 0 = \partial_z^2\tilde\mu_0, &\text{in }\Gamma_{b,\ell},\label{MB-eq:t1Inner0}
\end{align}
where~$\tilde\mu_0$ is defined in (\ref{MB-eq:tildeMu0}). This relation is satisfied by $\tilde u_0=U_b$, which is consistent
with our choice of initial data corresponding to the dressing of an admissible bilayer with the bilayer profile. Modulo this
form the next orders in the expansion take the form,
\begin{align}
    O(\eps^{-1}):&\quad \label{MB-eq:t1Inner1}0 = \partial_z^2\tilde\mu_1, &\text{in } \Gamma_{b,\ell},\\
    O(1):&\quad\label{MB-eq:t1Inner12}-V(s,S)\partial_z \tilde u_0  = \partial_z^2\tilde\mu_2+(\Hbl_0\partial_z+\Delta_S)\tilde\mu_1,&\text{in } \Gamma_{b,\ell},
\end{align}
where $\tilde\mu_1$ and~$\tilde\mu_2$ are defined in (\ref{MB-eq:tildeMu1}) and (\ref{MB-eq:tildeMu2}) respectively.
With the reduction $\tilde u_0=U_b$, the matching condition (\ref{MB-eq:MatchMu1}) reduces to the relations
$\partial_z\tilde\mu_1 \approx \partial_\nu\mu_0 = 0$ as~$z\longrightarrow\pm\infty$. Applying these to (\ref{MB-eq:t1Inner1}) we find that $\tilde\mu_1$ is independent of $z$,
i.e.,~$\tilde\mu_1 = \tilde \mu_1(s,S,\tau)$.
Similarly, we simplify the inner expression for $\tilde\mu_1$ in~(\ref{MB-eq:tildeMu1}) which yields the expression
\begin{equation}\label{MB-eq:u1tilde}
	\tilde u_1 =  \tilde \mu_1\Phi_{b,2}-\eta_dL_{b,0}^{-2}W'(U_b),
\end{equation}
where~$\Phi_{b,2}$ is defined in (\ref{e:Phi_bj-def}). Since $\tilde\mu_1$ is independent of $z$, the relation (\ref{MB-eq:t1Inner12}) reduces to
\begin{align}
\label{MB-eq:t1Inner2}
O(1):& \quad  -V(s,S)\partial_z U_b  = \partial_z^2\tilde\mu_2+\Delta_S \tilde\mu_1,& \text{in } \Gamma_{b,\ell}.
\end{align}

To obtain interfacial jump conditions for $\mu_1$ we introduce $\hat U_b:=U_b-b_->0$, and integrate (\ref{MB-eq:t1Inner2}) twice from $0$ to $z$ we
obtain the relation
\beq\label{MB-eq:t1innerMu2}
    \tilde\mu_2(z) = \tilde\mu_2(0) -V(s,S)\int_0^z\hat U_b(t)\,dt+z\left(\partial_z\tilde\mu_2(0)+{V(s,S)\hat U_b(0)}\right) +\frac{z^2}{2}\Delta_S\tilde\mu_1(s,S).
\eeq
Comparing to the jump condition (\ref{MB-eq:MatchMu2}), and recalling that $\mu_0=0$, we deduce that $\Delta_s \tilde\mu_1=0$. Since
$\Gamma_b$ is closed, this implies that $\tilde\mu_1=\tilde\mu_1(s,\tau)$ is constant in $S$. Using this information, we integrate (\ref{MB-eq:t1Inner2}) with
respect to $z$ over $\RR$.  As $U_b$ is homoclinic we obtain the relationship
\begin{align}
      \lim_{z\rightarrow\infty}\partial_z\tilde\mu_2(z)-\lim_{z\rightarrow-\infty}\partial_z\tilde\mu_2(z) = 0\label{MB-eq:t1KeyIdent2},
\end{align}
which when reported to the matching condition (\ref{MB-eq:MatchMu2}) yields the key outer interfacial jump relations
\begin{align}
    &\llbracket \mu_2\rrbracket = 0,\quad\llbracket \partial_\nu\mu_1\rrbracket = 0\label{MB-eq:mu1nJumpCond}.
\end{align}
Coupling these boundary conditions with the elliptic problem ({\ref{MB-eq:t1Outer1}) typically yields a Mullins-Sekerka problem for the long-range
transport of material; however the homogeneous jump conditions imply that $\Delta \mu_1=0$ in all of $\Omega$, which subject to
the exterior boundary conditions implies that the far-field chemical potential, $\mu_1$,  is spatially constant; however it remains a function of time.

To extract the normal velocity we return the reduction $\tilde u_0= U_b$ to (\ref{MB-eq:tildeMu2}), so that $P_0$ reduces to $L_{b,0}$
and the chemical potential takes the form
\begin{align}\label{MB-eq:t1tildeMu2exp}
	\tilde\mu_2  =&	L_{b,0}^2\tilde u_2-L_{b,0}\tilde Q_2
    +(\Hbl_0(s,S)\partial_z +\Delta_S-W'''(U_b)\tilde u_1+\eta_1)(L_{b,0} \tilde u_1+\Hbl_0(s,S)U_b')\\
    &+\eta_dW''(U_b)\tilde u_1,\nonumber
\end{align}
where we have introduced the quantity $\tilde Q_2:=Q_2-L_{b,0}\tu_2.$
Integrating (\ref{MB-eq:t1Inner2}) from~$z=-\infty$ to~$z=0$, using the matching condition~(\ref{MB-eq:MatchMu2}), and
recalling that $\mu_1$ is constant, $\Delta_S\tilde\mu_1=0$, and $\hat U_b\rightarrow 0$ as $z\rightarrow\pm\infty$ yields the expression
\beq\label{MB-eq:t1VhatUb}
    {V(s,S)\hat U_b(0)} = \lim_{z\rightarrow-\infty}\partial_z\tilde\mu_2-\partial_z\tilde\mu_2(0) = \partial_\nu\mu_1-\partial_z\tilde\mu_2(0)= -\partial_z\mu_2(0).
\eeq
Using~(\ref{MB-eq:t1VhatUb}) to replace~$V(s)\hat U_b(0)$ in equation~(\ref{MB-eq:t1innerMu2}) yields
\beq\label{MB-eq:t1tildeMu2}
    \tilde\mu_2(z) = \tilde\mu_2(0)-V(s,S)\int_0^z\hat U_b(t)\,dt.
\eeq
Replacing~$\tilde\mu_2$ in~(\ref{MB-eq:t1tildeMu2}) with its expression from~(\ref{MB-eq:t1tildeMu2exp}) and solving for~$L_{b,0}^2\tilde u_2$ yields
an expression for $\tilde u_2$
\begin{align}
\label{MB-eq:t1SolveCond}
    L_{b,0}^2\tilde u_2 = & L_{b,0}\tilde Q_1-\left(\Hbl_0\partial_z-W'''(U_b)\tilde u_1+\eta_1\right)\left(L_{b,0}\tilde u_1 + \Hbl_0\partial_z U_b\right)+ \nonumber \\
      &-\Delta_S H_0(s,S) \partial_z U_b- \eta_d W''(U_b)\tilde u_1+\tilde\mu_2(0)-V(s,S)\int_0^z\hat U_b(t)\,dt .
\end{align}
Fixing the values of $s$ and $S$, this equation has a solution~$\tilde u_2(s,S,\cdot) \in L^2(\RR)$ if and only if the right-hand side is
perpendicular to~$\ker L_{b,0}$, which is spanned by $\partial_z \phi_b.$ This solvability condition is enforced by selecting the value of $V(s,S)$. Since the terms
in (\ref{MB-eq:t1SolveCond}) are either functions of $z$ or of $s$ and $S$, we factor out the functions of $s$ and $S$, replace $U_b$ with $\phi_b$, and
take the inner product of (\ref{MB-eq:t1SolveCond}) with $\partial_z \phi_b $ in $L^2(\RR).$ Recalling that $\tilde u_1$, defined in~(\ref{MB-eq:u1tilde}), is even in $z$ and
the operator~$L_{b,0}$ preserves symmetry, parity considerations reduce the solvability condition to
\beq\label{MB-eq:t1SolvCondInnerProduct}
\Hbl_0  \Bigl( (L_{b,0}\tu_1,\partial_z^2 \phi_b)_{L^2(\RR)}+ (W'''(\phi_b) \tu_1, \partial_z \phi_b) - \eta_1\|\phi_b'\|_{L^2}^2\Bigr) +V \|\hat\phi_b\|_{L^2}^2 - \Delta_S \Hbl_0 \|\phi_b^\prime\|_{L^2}^2 =0.
\eeq
From (\ref{e:BLcp}) it is easy to verify that
\begin{align}
L_{b,0} \left(\frac{z}{2} \phi_b^\prime\right)=  \phi_b^{\prime\prime},\\
L_{b,0} \phi_b^{\prime\prime}=W'''(\phi_b)(\phi_b^\prime)^2.
\end{align}
Using these relations and the form (\ref{MB-eq:u1tilde}) of $\tu_1$, the coefficient of $\Hbl_0$ in (\ref{MB-eq:t1SolvCondInnerProduct}) reduces to
\begin{align}
(L_{b,0}\tu_1,\partial_z^2 \phi_b)_{L^2(\RR)}+ (W'''(\phi_b) \tu_1, \partial_z \phi_b) -\eta_1 \|\phi_b'\|_{L^2}^2 &=
 (L_{b,0}^2\tu_1, \frac{z}{2}\phi_b^\prime)_{L^2}  -\eta_1 \|\phi_b'\|_{L^2}^2,\nonumber \\
 &= \mu_1 m_b+\frac{1}{2}\left(\eta_1+\eta_2\right)\sigma_b.\nonumber
 \end{align}
Returning this reduction to (\ref{MB-eq:t1SolvCondInnerProduct}) and solving for the normal velocity we find
\beq
 V(s,S) = \frac{ \mu_1 m_b +\frac12(\eta_1+\eta_2)\sigma_b}{B_1} \Hbl_0 + \frac{\sigma_b }{B_1}\Delta_S\Hbl_0,
 \eeq
where here and above we have introduced the positive constants
\beq
    \label{e:MB-const} m_b:=\int_\mbbR \hat\phi_b\,dz>0, \quad
                                  B_1:=\norm{\hat \phi_b}_{L^2(\RR)}^2,\quad
                                 \sigma_b:=\norm{\phi_b'}_{L^2(\RR)}^2.
 \eeq
The sign of the coefficient of $\Hbl_0$ is  indeterminate, as $\eta_2$ can be negative and moreover the bulk chemical potential
$\mu_1$ varies temporally. To emphasize this fact we introducing the constants
\beq\label{e:nub-def2}
\mu_b^*=-\frac{1}{2}\left(\eta_1+\eta_2\right)\frac{\sigma_b}{m_b}, \hspace{0.5in}  \nu_b:= \frac{m_b}{B_1}, \hspace{0.4in} k_b:=\frac{\sigma_b}{B_1},
\eeq
and return the $S$ variable to its original scaling, obtaining the regularized curvature driven flow
\beq
V(s) = \nu_b \left(\mu_1-\mu_b^*\right)\Hbl_0 +\eps k_b \Delta_s \Hbl_0.\label{e:Vb-eq}
\eeq

To close the system and fully determine the normal velocity we evoke conservation of mass to specify the temporally varying value of
the bulk external chemical potential, $\mu_1$. The mass balance is determined by the interplay between the length of the
interface $\Gamma_b$ and the total mass of amphiphilic material. From
the form of $\tu_0$ and $\tu_1$ in (\ref{MB-eq:u1tilde}), we have the composite formulation
\beq\label{MB-eq:t1InnerSol}
    u(x,t) = U_b + \eps(\mu_1\Phi_{b,2}-\eta_dL_{b,0}^{-2}W'(U_b))+O(\eps^2),
\eeq
which has the spatially constant far-field asymptotic value
\beq\label{MB-eq:t1OuterSol}
    u(x,t) = b_-+\eps\frac{\mu_1}{\alpha_-^2}+O(\eps^2)\quad\text{in }\Omega\backslash\Gamma_{b,\ell},
\eeq
where~$\alpha_- := W''(b_-)>0$.
The gradient flow (\ref{e:FCH-eq}) conserves the total mass,
\beq
    M:=\int_\Omega u(x,t)-b_-\,dx = \int_\Omega u(x,0)-b_-\,dx.
\eeq
Using the form of the composite solution, (\ref{MB-eq:t1InnerSol}), we evaluate the integral over the reach and its complement,
\beq
    M = \eps\int_{\tilde\Gamma_{b,\ell}}\frac{\mu_1}{\alpha_-^2}\,dx+\int_{\Gamma_{b,\ell}}\hat U_b+ \eps(\mu_1\Phi_{b,2}-\eta_dL_{b,0}^{-2}W'(U_b))\,dx+O(\eps^2).
\eeq
Since $\Gamma_b$ is admissible, its area $|\Gamma_b|\sim O(1)$.  Changing to whiskered coordinates in the localized integral yields
\beq\label{MB-eq:M}
    M = \eps\left(|\Omega|\frac{\mu_1}{\alpha_-^2}+\int_{\Gamma_b}\intreach \hat \phi_b\,dz\,ds\right)+O(\eps^2).
\eeq
Our choice of initial data implies that the mass can be rescaled as $M=\eps \hat M+O(\eps^2)$. We also expand the surface area
\beq
\label{MB-eq:GammaExp}
|\Gamma_b| = \gamma_{b,0}+\eps\gamma_{b,1}+O(\eps^2).
\eeq
Evaluating the integrals in equation (\ref{MB-eq:M}) and solving for $\mu_1$ yields the expression
\beq\label{MB-eq:Mu1Gamma0}
    \mu_1 = \frac{\alpha_-^2}{|\Omega|}\left(\hat M-\gamma_{b,0}m_b\right) + O(\eps),
\eeq
where~$m_b$ is defined in~(\ref{e:MB-const}).
On the other hand, the area of a smooth curve subject to normal velocity $V$ evolves according to
\beq\label{MB-eq:gammaGrowthBl}
	\frac{\partial|\Gamma_b|}{\partial \tau}=\int_{\Gamma_b} V(s)\Hbl_0(s)\,ds,
\eeq
and for the normal velocity (\ref{e:Vb-eq}) this reduces to the result presented in (\ref{e:Gb-grow}).
Taking the time derivative  of (\ref{MB-eq:Mu1Gamma0}), using (\ref{MB-eq:gammaGrowthBl}) to eliminate $\frac{d}{d\tau}\gamma_{b,0}$,
the normal velocity (\ref{e:Vb-eq}) drives the bulk chemical potential according to
\beq
\frac{d\mu_1}{d\tau} =- \frac{\alpha_-^2m_b}{|\Omega|}\left(
\nu_b(\mu_1-\mu_b^*)\int_{\Gamma_b} \Hbl_0^2(s)\,ds -\eps k_b \int_{\Gamma_b}|\nabla_s \Hbl_0|^2\, ds\right)+O(\eps^2).
  \label{e-mu1b}
\eeq

The coupled system (\ref{e:Vb-eq}) and (\ref{e-mu1b}) prescribes the interfacial evolution for the dressing of an admissible codimension one interface with
a shifted bulk value of $u$ given by (\ref{MB-eq:t1OuterSol}).

\section{Geometric evolution of codimension two structures}\label{Sec:Filament}
In this section we derive the geometric evolution of admissible codimension two curves, called filaments, in $\RR^3$ under the~$H^{-1}$ gradient flow (\ref{e:FCH-eq}).
As remarked in the introduction, we make the assumption that the combined length of all the filament curves scales as $O(\eps^{-1})$ so that the combined mass of the codimension two
structures is $O(\eps)$. This gives a comparable masses to filament and bilayer structures, so they may contribute to the mass balance at the same order of magnitude.
Codimension two structures are much less studied than codimension one structures, however our analysis leads to a qualitatively similar result: a surface diffusion regularized
curvature-vector driven normal flow that may be curve lengthening or curve shortening depending upon the value of the spatially constant far-field chemical potential.

\subsection{Admissible codimension two curves and their dressings}

Given a smooth, closed, non-self intersecting one-dimensional manifold $\Gamma_\fil$ immersed in $\Omega\subset\RR^3$,
and parameterized by the map $s\in \cS_\fil\mapsto\zeta_\fil(s)\in\Omega$, we may uniquely decompose points $x$ near $\Gamma_\fil$ as
\beq
\label{e:cd2-wk}
  x = \rho_\fil(s,z_1,z_2) = \zeta_\fil(s)+\eps\left(z_1N_1(s)+z_2N_2(s)\right),
\eeq
where $N_1(s)$ and $N_2(s)$ are orthogonal unit vectors which are also orthogonal to the tangent vector $\zeta_\fil^\prime(s)$, defined by
\beq
    \frac{\partial \bN^i}{\partial s} = -\kappa_i\bT,\quad i=1,2,
\eeq
where
\beq\label{CS-eq:VecKappa}
    \vec{\kappa}(s,t) := (\kappa_1,\kappa_2)^t,
\eeq
is the normal curvature vector with respect to~$\{\bN^1,\bN^2\}$.
\begin{definition}
\label{def:admissible_f}
 For any $K,\ell > 0$ the family, $\cG^\fil_{K,\ell}$, of admissible curves is comprised of closed (compact and without boundary), oriented 1 dimensional curves $\Gamma_\fil$ embedded in $\RR^3$, which are far from self intersection and with a smooth second fundamental form. More precisely,
\vskip 0.01in
\begin{tabular}{lp{5.0in}}
 (i)& The $W^{4,\infty}(\cS_\fil)$ norm of the 2nd Fundamental form of $\Gamma_\fil$ and its principal curvatures
are bounded by $K$.\\
 (ii)& The whiskers of length $3\ell < 1/K$, in the unscaled distance, defined for each $s_0\in \cS_\fil$
by, $w_{s_0} := \{x : s(x) = s_0,~|\vz(x)|<3\ell/\eps\}$, neither intersect each-other nor $\partial\Omega$ (except when considering periodic boundary conditions).\\
(iii)& The length, $|\Gamma_\fil|$, of $\Gamma_\fil$ is bounded by $K/\eps$.
\end{tabular}
\end{definition}

For an admissible codimension two curve $\Gamma_\fil$ the change of variables $x \rightarrow \rho(s,\vz)$ given by (\ref{e:cd2-wk}) is a $C^4$ diffeomorphism
on the reach of $\Gamma_\fil$, defined as the set
\beq
\label{e:whiskNbd-f}
\Gamma_{\fil,\ell} := \left\{ \rho(s,\vz) \in \mbbR^3 \Bigl| s \in \cS, -\ell/\eps \le |\vz| \le \ell/\eps \right\}\subset\Omega.
\eeq
where $\vz:=(z_1,z_2)$. We introduce $R(x)=|\vz(x)|$ which denotes the scaled distance of $x$ to $\Gamma_\fil.$ Within the reach the
cartesian Laplacian admits the local form
\beq\label{e:LaplacianPr}
    \eps^2\Delta_x = \Delta_{\vz}-\eps\vec{\kappa}\cdot\nabla_{\vz}+\eps^2(\partial_s^2-(\vz\cdot\vec{\kappa})\vec{\kappa}\cdot\nabla_{\vz})+O(\eps^3),
\eeq
where the lower order terms are immaterial for the analysis. The Jacobian of the change of variables (\ref{e:cd2-wk}) takes the form
\beq\label{e:JacobianFil}
    J = \eps^2-\eps^2\vz\cdot\vec{\kappa}.
\eeq
If the underlying curve $\Gamma_\fil$ evolves in time, then its normal velocity vector~$\vV = (V_1,V_2)$ of~$\Gamma_\fil$ at a point~$s(t)$ takes the form
 \begin{align}\label{e:VDefFil}
    &V_1:= -\eps\frac{\partial z_1}{\partial t}+\eps z_2N^2\cdot\frac{\partial N^1}{\partial t},\\
    &V_2:= -\eps\frac{\partial z_2}{\partial t}+\eps z_1N^1\cdot\frac{\partial N^2}{\partial t}.
 \end{align}
 The terms $z_2N^2\cdot\frac{\partial N^1}{\partial t}$ and $z_1N^1\cdot\frac{\partial N^2}{\partial t}$ reflect lower-order contributions to the normal velocity
 induced by the rotational motion of the normal vectors to $\Gamma_\fil(t)$.
See \cite{Dai2015competitive} for further details.

\begin{definition} Given an admissible codimension two curve $\Gamma_\fil\in\cG_{K,\ell}^\fil$ and a smooth function~$f:\RR_+\rightarrow\RR$ which tends to a constant
value~$f_\infty$ at an $O(1)$ exponential rate as~$R\rightarrow\infty$, we define $f_{\Gamma_\fil}\in H^2(\Omega)$, called the {\bf dressing of $\Gamma_\fil$ with~$f$}, according to the rule
\beq
    f_{\Gamma_\fil}(x):=f(R(x))\chi(R(x)/\ell)+f_\infty(1-\chi(R(x)/\ell)),
\eeq
where~$\chi:\RR\rightarrow\RR$ is a fixed, smooth cut-off function taking values one on~$[0, 1]$, and zero on $[2,\infty)$.
By abuse of notation we will drop the $\Gamma_\fil$ subscript when doing so creates no confusion.
\end{definition}

Within the reach $\Gamma_{\fil,\ell}$ the Cartesian Laplacian reduces formally at leading order to the two-dimensional Laplacian in $\vz$, which may be written
in turn in polar coordinates in $R$.  We may eliminate the dominant terms in quadratic component of the (\ref{e:FCHFE}) by taking $u$ at leading order
to be the dressing of the codimension two profile~$\phi_\fil(R)$, defined as the solution of
\beq
\label{e:FLcp}
  \partial_R^2 \phi_\fil+\frac{1}{R}\partial_R\phi_\fil - W'(\phi_\fil)=0,
\eeq
subject to the boundary conditions $\partial_R \phi_\fil(0)=0$ and $\phi_\fil\to b_-$ as $R\to\infty$.
We denote  the dressing of $\Gamma_\fil$ with $\phi_\fil$ by $U_\fil$. As in the codimension one case, we introduce the operator
\beq \label{e:Lf0-def}
L_{\fil,0} := \partial_R^2 +\frac{1}{R}\partial_R-W''(\phi_\fil),
\eeq
corresponding to the linearization of (\ref{e:FLcp}) about $\phi_{\fil}$.
The operator $L_{\fil,0}$ is the radially symmetric reduction of the associated cylindrical Laplacian,
\beq\label{CS-eq:L}
    L_\fil:=\partial_R^2+\frac{1}{R}\partial_R+\frac{1}{R^2}\partial_\theta^2-W''(U_\fil).
\eeq
This operator is self-adjoint in the usual $R$-weighted $L^2(\RR_+)$ inner product,
\beq \langle f, g \rangle_R:= \int_0^\infty f(R)g(R) R\,dR.
\eeq
Moreover, the translational eigenfunctions $\{\phi_\fil^\prime(R)\cos(\theta), \phi_\fil^\prime(R)\sin(\theta)\}$ lie in the kernel of $L_\fil$ and their associated dressings of $\Gamma_\fil$ agree with $\{\partial_{z_1}U_\fil, \partial_{z_2}U_\fil\}$ respectively, up to exponentially small terms.
For each $m\in\NN$, we define the spaces
\beq\label{CS-eq:ZmSpace}
    \cZ_m:=\{f(R)\cos(m\theta)+g(R)\sin(m\theta)~\big|~ f,g\in C^\infty(0,\infty),m\in\mathbb{N}\}.
\eeq
These spaces are invariant under the operator~$L_\fil$, and mutually orthogonal in~$L^2(\Omega)$. Moreover, on these spaces~$L_\fil$ reduces to
\beq
    L_\fil(f(R)\cos(m\theta)+g(R)\sin(m\theta)) = \cos(m\theta)L_{\fil,m}f+\sin(m\theta)L_{\fil,m}g,
\eeq
where
\beq\label{CS-eq:Lm}
    L_{\fil,m}:=\frac{\partial^2}{\partial R^2}+\frac{1}{R}\frac{\partial }{\partial R}-\frac{m^2}{R^2}-W''(U_\fil).
\eeq
Each operator~$L_{\fil,m}$ is self-adjoint in the~$R$-weighted inner product, and
the operator~$L_{\fil,1}$ has a 1-dimensional kernel spanned by its ground state $\partial_R \phi_\fil>0$. For~$m > 1$ we observe that~$\langle L_{\fil,m}f,f\rangle_{R}<\langle L_{\fil,1}f,f\rangle_{R}$ and since $L_{\fil,1}\leq0$ we deduce that $L_{\fil,m}<0$ for $m>1$, and is boundedly invertible.
We denote the eigenfunctions and eigenvalues of $L_{\fil,m}$ by $\{\psi_{\fil,m,j}\}_{j=0}^\infty$ and $\{\lambda_{\fil,m,j}\}_{j=0}^\infty$, respectively.
We address the kernel of $L_{\fil,0}$ with the following assumption.
\begin{asmp}
\label{A-Lf}
We assume that the operator $L_{\fil,0}$ has no kernel and a one-dimensional positive eigenspace spanned by $\psi_{\fil,0,0}$.
\end{asmp}
With this assumption we may define the functions
\beq
\Phi_{\fil,j}:=L_{\fil,0}^{-j} 1, \label{e:Phi_fj-def}
\eeq
for $j=1, 2$ and their $\Gamma_{\fil}$ dressings, also denoted $\Phi_{\fil,1}$ and $\Phi_{\fil,2}$.

\subsection{Inner and outer expansions}
Considering initial data that is close to a filament dressing of an admissible curve, $\Gamma_\fil(0)\in\cG_{K,\ell}^\fil$, embedded in $\Omega\subset\RR^3$.
In the far-field, $\tilde\Gamma_{\fil,\ell}$, the outer solution $u$ has the expansion
\begin{align}
    \label{e:uOuterSolFil}u(x,t) &= u_0(x,t)+\eps u_1(x,t)+\eps^2 u_2(x,t)+O(\eps^3),
    \end{align}
and within the reach $\Gamma_{\fil,\ell}$, we incorporate a two-scale tangential expansion, introducing the variable $ S=\frac{s}{\sqrt{\eps}}$, and the inner spatial expansion takes the form
\begin{align}
    \label{e:uInnerSolFil}u(x,t)&=\tilde u(s,S,\vz,\tau) = \tilde u_0(s,S,\vz,\tau)+\eps\tilde u_1(s,S,\vz,\tau)+\eps^2 \tilde u_2(s,S,\vz,\tau)+O(\eps^3).\
\end{align}
The time derivative of the inner density function~$\tilde u$, defined in~(\ref{e:uOuterSolFil}), combined with the normal velocity,~(\ref{e:VDefFil}), takes the form
\beq\label{MB-eq:ut}
    \frac{\partial \tilde u}{\partial t} = -\eps^{-1}\bV\cdot\nabla_\vz\tilde u+\frac{\partial \tilde u}{\partial \tau}\frac{\partial \tau}{\partial t}.
\eeq

For a whisker identified by $s\in \cS_\fil$, with base point $x=\rho_\fil(s,0)\in\Gamma_\fil$  we choose vectors $\bn,\bm\in \text{span}\{\bN^1,\bN^2\}$ in the normal plane of
$\Gamma_\fil$ at $x$, and choose $\theta$ so that
\beq
    \bn = \cos(\theta)\bN^1+\sin(\theta)\bN^2.
\eeq
The usual directional derivative along~$\bn$ is denoted
\beq\label{MP-eq:DirectionalDerivative}
    \partial_\bn := \bn\cdot\nabla_x = \cos(\theta)\bN^1\cdot\nabla_x+\sin(\theta)\bN^2\cdot\nabla_x,
\eeq
and for~$f\in C^\infty(\Omega/\Gamma_\fil)$ we introduce the~$\bn,\bm$ limit
\beq
    \partial_\bn^jf^\bm(x):=\lim_{h\rightarrow0^+}(\bn\cdot\nabla_x)^jf(x+h\bm,t)\quad\text{for all }j\geq0,
\eeq
and the limit of the gradient
\beq
    \nabla_xf^\bm(x):=\lim_{h\rightarrow0^+}\nabla_xf(x+h\bm,t),
\eeq
where the limit exists. If~$f\in C^1(\Omega)$ then the normal derivative of~$f$ will satisfy
\beq
    \partial_\bn f^{-\bm} = \partial_\bn f^{\bm}.
\eeq
This motivates the following definition of the jump condition.
\begin{definition}\label{MP-def:JumpCond}
Given a radial function~$f:=f(R)$ localized on~$\Gamma_\fil$, we define {\em the jump of~$f$ across a given whisker} by
\beq\label{MP-eq:JumpCondDef}
    \llbracket \partial_\bn f^{\bm}\rrbracket_{\Gamma_\fil} (x) := \partial_\bn f^{\bm}(x) - \partial_\bn f^{-\bm}(x)
\eeq
which is zero when~$f$ has a smooth extension through~$\Gamma_\fil$.
\end{definition}

With this notation we develop matching conditions
\beq
 \lim\limits_{R\to 0^+} u(x+\eps R\bn,t) = \lim\limits_{R\to\infty} \tu(s,S,R,\theta,t).
 \label{MP-eq:Matching}
 \eeq
Expanding the left-hand side yields the following expression
\begin{align}\label{MP-eq:OuterMatch}
   u(x+\eps R\bn)=& u_0^{\bn}(x,t)+\eps\left(u_1^{\bn}(x,t)+R\partial_\bn u_0^{\bn}(x,t)\right)\\
   &+\eps^2\left(u_2^{\bn}(x,t)+R\partial_\bn u_1^{\bn}(x,t)+\frac{1}{2}R^2\partial_\bn^2 u_0^{\bn}(x,t)\right)+ O(\eps^3),\nonumber
\end{align}
where $u_i^{\bn}$ denotes the limit of the left-hand side of (\ref{MP-eq:Matching}) as $\eps R\to0^+$.
Equating orders of~$\eps$ for the matching condition~(\ref{MP-eq:Matching}) yields
\begin{align}
    u_0^{\bn} &= \lim_{R\rightarrow\infty}\tilde u_0(s,S,R,\theta,\tau),\label{MP-eq:MatchU0}\\
    u_1^{\bn}+R\partial_\bn u_0^{\bn} &= \lim_{R\rightarrow\infty}\tilde u_1(s,S,R,\theta,\tau).
\end{align}

The chemical potential, defined in (\ref{e:mu}) admits similar inner and outer expansions. The outer expansion is identical to that
for the codimension one case, see (\ref{MB-eq:Mu0}) and (\ref{MB-eq:Mu1}). To obtain the inner expression for the chemical potential we first note that in the multiscale tangential variables the Laplacian expression (\ref{e:LaplacianPr}) takes the form
\beq
    \eps^2\Delta_x=
    \Delta_{\vz}-\eps\left(\vec{\kappa}\cdot\nabla_{\vz}+\partial_S^2\right)+\eps^2(\partial_s^2-(\vz\cdot\vec{\kappa})\vec{\kappa}\cdot\nabla_{\vz}+{2\vz\cdot\vec{\kappa}\partial_S^2})+O(\eps^3),
\eeq
Introducing the nonlinear opeartors $P$ and $Q$, the inner chemical potential is written as
\beq\label{MP-eq:InnerMu}
    \tilde{\mu}(x,t) = P(\tu)Q(\tu)+\eps\eta_dW'(\tu_0)
\eeq
where $ P$ admits the expansion $P=P_0+\eps  P_1+\eps^2 P_2+ \ldots$, with
\begin{align}
\label{MP-e:P0}
P_0 = & -\Delta_\vz + W''(\tu_0), \\
\label{MP-e:P1}
 P_1 = & \vec{\kappa}\cdot \nabla_\vz +\partial_S^2+ W'''(\tu_0)\tu_1 - \eta_1, \\
\label{MP-e:P2}
 P_2 = & -(\vz\cdot\vec{\kappa})\vec{\kappa}\cdot\nabla_\vz + {2\vz\cdot\vec{\kappa}\partial_S^2} - \partial^2_s + W'''(\tu_0)\tu_2 +\frac12 W^{(4)}(\tu_0)\tu_1^2,
\end{align}
and simliarly $Q=Q_0+\eps Q_1+\eps^2 Q_2 +\ldots$, where
\begin{align}
\label{MP-e:Q0} Q_0 = & -\Delta_\vz\tu_0 + W'(\tu_0), \\
\label{MP-e:Q1}  Q_1 = & \vec{\kappa}\cdot\nabla_\vz \tu_0 +\partial_S^2\tu_0 +\left(-\Delta_\vz+ W''(\tu_0)\right)\tu_1, \\
\label{MP-e:Q2} Q_2 = & -(\vz\cdot\vec{\kappa})\vec{\kappa}\cdot\nabla_\vz\tu_0 - {2\vz\cdot\vec{\kappa}\partial_S^2}\tu_0 - \partial^2_s\tu_0 +(\vec{\kappa}\cdot\nabla_\vz+\partial^2_S)\tu_1
+ (-\Delta_\vz+W''(\tu_0))\tu_2\\
& -\frac12 W'''(\tu_0)\tu_1^2.\nonumber
\end{align}
With these reductions we expand the inner chemical potential as
\begin{align}
\label{MP-eq:tildeMu0}\tilde\mu_0=& P_0Q_0\\
\label{MP-eq:tildeMu1}\tilde\mu_1 =& P_1Q_0+P_0Q_1 +\eta_d W'(\tu_0), \\
\label{MP-eq:tildeMu2}\tilde\mu_2 =& P_0Q_2+P_1Q_1 +P_2Q_0+\eta_d W''(\tu_0)\tu_1 .
\end{align}
The relevant matching conditions for the chemical potential extend to second order in $\eps$:
\begin{align}
    \mu_0^{\bn}(x,t) &= \lim_{R\rightarrow\infty}\tilde \mu_0(s,S,R,\theta,t),\label{MP-eq:MatchMu0}\\
    \mu_1^{\bn}(x,t)+R\partial_\bn \mu_0^{\pm}(x,t) &= \lim_{R\rightarrow\infty}\tilde \mu_1(s,S,R,\theta,t),\label{MP-eq:MatchMu1}\\
    \mu_2^{\bn}(x,t)+R\partial_\bn \mu_1^{\bn}(x,t)+\frac{1}{2}R^2\partial_\bn^2\mu_0^\bn(x,t) &= \lim_{R\rightarrow\infty}\tilde \mu_2(s,S,R,\theta,t),\label{MP-eq:MatchMu2}
\end{align}

\subsection{Time scale $\tau=\eps t$ : quenched vector-curvature driven flow}
The analysis of the outer expansion of the chemical potential is identical to the codimension one case, and we find at leading order that
$u_0=b_-$, $\mu_0=0$, while at $O(\eps)$ we obtain
 \beq\label{MP-eq:t1Outer1}
 	\Delta_x\mu_1=0\quad\text{in }\Omega\backslash \Gamma_{\fil}.
 \eeq
In the inner region we supplement the inner expansions (\ref{e:uInnerSolFil}) and (\ref{MP-eq:InnerMu}) with the inner expression of the Laplacian (\ref{e:LaplacianPr}). At leading order in $\eps$ we find
\begin{align}
    O(\eps^{-2}): &\quad 0 = \Delta_\vz\tilde\mu_0,\quad\text{in } \Gamma_{\fil,\ell} \label{MP-eq:t1Inner0},
 \end{align}
 where $\tilde\mu_0$ is defined in (\ref{MP-eq:tildeMu0}). This equation is consistent with the choice of initial data $\tilde u_0=\phi_\fil$ which implies that $\tilde\mu_0=0$  via the matching conditions.  With this reduction the subsequent orders become
 \begin{align}
   O(\eps^{-1}):&\quad 0=\Delta_\vz \tilde\mu_1,&\text{in }\Gamma_{\fil,\ell}, \label{MP-eq:t1Inner1}\\
    O(1):&\quad -\bV\cdot\nabla_\vz\tilde u_0 = \Delta_\vz\tilde\mu_2-\left(\vec\kappa\cdot\nabla_\vz+\partial_S^2\right)\tilde\mu_1, &\text{in }\Gamma_{\fil,\ell}. \label{MP-eq:t1Inner12}
\end{align}
where $\tilde\mu_1$ and $\tilde\mu_2$ are defined in equations (\ref{MP-eq:tildeMu1}) and (\ref{MP-eq:tildeMu2}), respectively.
The combined system (\ref{MP-eq:t1Outer1}) and (\ref{MP-eq:t1Inner1}) couples through the matching condition
(\ref{MP-eq:MatchMu1}). Since $\mu_0=0$ we deduce from (\ref{MP-eq:MatchMu1})
that $\tilde\mu_1$ is bounded as $R\to\infty$, and hence from (\ref{MP-eq:t1Inner1}) that $\tilde\mu_1$
is constant in $\vz$. In particular $\tilde\mu_1=\tilde\mu_1(s,S,\tau)\approx \mu_1$. Since $\tilde u_0=U_\fil$, equation (\ref{MP-eq:tildeMu1}) for $\tilde\mu_1$
reduces to a linear equation for $\tu_1$,
\beq\label{MP-eq:t1tildeMu1}
    L_\fil^2\tilde u_1=\tilde\mu_1 -\eta_dW'(U_\fil).
\eeq
By Assumption\,\ref{A-Lf} we know that $\ker L_\fil\subset\ker L_{\fil,1}\subset\cZ_1$, defined in (\ref{CS-eq:ZmSpace}), while the right-hand side of (\ref{MP-eq:t1tildeMu1}) lies in  $\cZ_0$.
Since the spaces~$\cZ_m$ are mutually orthogonal, we may solve for $\tilde u_1$,
\beq\label{MP-eq:TildeU1}
	\tilde u_1 = \mu_1\Phi_{\fil,2}-\eta_dL_\fil^{-2}W'(U_\fil),
\eeq
where $\mu_1$ is a spatial constant and~$\Phi_{\fil,2}$ is defined in (\ref{e:Phi_fj-def}).
With this simplification equation (\ref{MP-eq:t1Inner12}) becomes
\beq\label{MP-eq:t1Inner2}
    O(1): -\bV\cdot\nabla_\vz U_\fil = \Delta_\vz\tilde\mu_2- \partial_S^2\tilde\mu_1.
\eeq
To impose interfacial matching conditions for $\tilde\mu_2$ we solve (\ref{MP-eq:t1Inner2}) by expanding $\tilde\mu_2$ in $(R,\theta)$ inner-polar coordinates
associated to the spaces $\{\cZ_m\}_{m=0}^\infty$ as
\beq\label{MP-eq:t1TildeMu2Full}
    \tilde\mu_2 =  A_1(s,S,R)\cos\theta + B_1(s,S,R)\sin\theta+\bar C(s,S,R)+\xi(s,S,R,\theta),
\eeq
where
\beq
    \xi(s,S,R,\theta) := \sum_{m=2}^\infty\left(A_m(s,S,R)\cos(m\theta)+ B_m(s,S,R)\sin(m\theta)\right).
\eeq
We observe that
\beq
\bV\cdot\nabla_\vz U_\fil = \partial_R U_\fil(R)\left(V_1 \cos\theta +V_2 \sin\theta\right) \in \cZ_1,
\eeq
while $\partial_S^2 \tilde\mu_1\in\cZ_0.$  We  project (\ref{MP-eq:t1Inner2}) onto $\cZ_m$ where $\Delta_\vz= \partial_R^2 +\frac{1}{R}\partial_R - \frac{m^2}{R^2}$ and arrive at the system
\begin{align}
    \partial_R^2 C+\frac{1}{R}\partial_RC &= \partial_S^2\tilde\mu_1(s,S),\label{MP-eq:t12C}\\
    \partial_R^2 A_1+\frac{1}{R}\partial_R A_1+\frac{1}{R^2}A_1 &= V_1(s,S)\partial_R U_\fil(R),\label{MP-eq:t12A1}\\
    \partial_R^2 B_1+\frac{1}{R}\partial_R B_1+\frac{1}{R^2}B_1 &= V_2(s,S)\partial_RU_\fil(R),\label{MP-eq:t12B1}
  \end{align}
 plus homogeneous equations for $\{A_m, B_m\}_{m=2}^\infty$ that have non-singular solutions $A_m=a_m(s,S)R^m$ and $B_m=b_m(s,S)R^m$. The $\cZ_0$
 equation has solution
 \beq \label{MP-eq:Z0}
 C(s,S,R)=C_0(s,S)+\frac{R^2}{4}\partial_S^2 \tilde\mu_1(s,S),
 \eeq
while the non-homogeneous equations, (\ref{MP-eq:t12A1}) and~(\ref{MP-eq:t12B1}), have the solutions
\begin{align}
    A_1(s,S,R) &= a_{1}(s,S)R-a(R)V_1(s,S),\label{MP-eq:t12ASol}\\
    B_1(s,S,R) &= b_{1}(s,S)R-a(R)V_2(s,S),\label{MP-eq:t12BSol}
\end{align}
where $a(R)$ is the solution of the non-homogeneous ordinary differential equation
\beq\label{MP-eq:t12a}
    a''+\frac{1}{R}a'-\frac{1}{R^2}a = \partial_RU_\fil(R),
\eeq
which enjoys the explicit formula
\beq\label{MP-eq:aR}
    a(R) = \frac{1}{R}\int_0^R r\hat U_\fil(r)\,dr,
\eeq
where we have introduced $\hat U_\fil:= U_\fil-b_-.$ In particular, $a(R)\to 0$ as $R\to\infty.$

From the matching condition~(\ref{MP-eq:MatchMu2}) we see that~$\tilde\mu_2$ grows at most linearly as~$R\rightarrow\infty$ and
\beq\label{MP-eq:mu2Growth}
    \lim_{R\rightarrow\infty}\frac{\partial\tilde\mu_2}{\partial R} = \partial_\bn \mu_1^{\bn}=\cos\theta\bN^1\cdot\nabla_x\mu_1^{\bn}+\sin\theta\bN^2\cdot\nabla_x\mu_1^{\bn},
\eeq
where the second equality follows from the definition of the directional derivative along~$\bn$, given in~(\ref{MP-eq:DirectionalDerivative}).
Taking the~$R$ derivative of~(\ref{MP-eq:t1TildeMu2Full}) and using the results above yields
\beq\label{MP-eq:t1tildeMu2dR}
    \frac{\partial \tilde\mu_2}{\partial R} = \frac{R}{2} \partial_S^2\tilde\mu_1 +(a_1-a'(R)V_1)\cos\theta+(a_2-a'(R)V_2)\sin\theta+\frac{\partial\xi}{\partial R}.
\eeq
Projecting each of (\ref{MP-eq:t1tildeMu2dR}) and (\ref{MP-eq:mu2Growth}) onto $\cZ_m$ and matching terms, we conclude that $\xi=0$ and
$\partial_S^2\tilde\mu_1=0$, and in particular
\beq\label{MP-eq:t1tildeMu2Parity}
    \frac{\partial \tilde\mu_2}{\partial R}(s,S,R,\theta,\tau) = -\frac{\partial \tilde\mu_2}{\partial R}(s,S,R,\theta+\pi,\tau).
\eeq
This latter result, substituted into (\ref{MP-eq:mu2Growth}) yields the no-jump condition across the curve $\Gamma_\fil$
\beq\label{MP-eq:t1JumpCond}
    \llbracket \partial_\bn \mu_1^{\bn}\rrbracket_{\Gamma_\fil} = 0,
\eeq
for any choice of normal vector~$\bn$. As the codimension two curve $\Gamma_\fil$ has zero capacity, it follows from the
zero-jump condition and (\ref{MP-eq:t1Outer1}) that $\mu_1$ has a harmonic extension to all of $\Omega$, and hence is spatially constant.
In particular $\nabla_x\mu_1^{\bn}=0$ for all choices of direction $\bn.$
With these reductions $\tilde\mu_2$ takes the form
\beq\label{MP-eq:t0TildeMu1}
    \tilde\mu_2 = C_0(s,S)-a(R)\left(V_1(s,S)\cos\theta+V_2(s,S)\sin\theta\right).
\eeq


To determine the normal velocity we substitute $\tilde u_0=U_\fil$  into the expression (\ref{MP-eq:tildeMu2}) for $\tilde\mu_2$, so that $P_0$ reduces to $L_\fil$ and the chemical potential takes the form
\begin{align}\label{MP-eq:t1Mu2U2}
	\tilde \mu_2 =& L_\fil^2\tilde u_2 -L_\fil\tilde Q_2 +(\vec{\kappa}\cdot \nabla_\vz +\partial_S^2+ W'''(\tu_0)\tu_1 - \eta_1)(-L_\fil\tilde u_1+\vec{\kappa}\cdot\nabla_\vz U_\fil)\\
&+
	\eta_dW''(U_\fil)\tilde u_1,\nonumber
\end{align}
where we have introduced $\tilde Q_2:= Q_2-L_\fil \tilde u_2$.
To solve equation~(\ref{MP-eq:t1Mu2U2}) for~$\tilde u_2$ we rewrite it in the  form
\beq\label{MP-eq:t1Tildeu2}
    L_\fil^2\tilde u_2 = \tilde\mu_2-\cQ+L_\fil\tilde Q_2,
\eeq
where
\beq\label{MP-eq:t1Q}
    \cQ:=(\vec{\kappa}\cdot \nabla_\vz +\partial_S^2+ W'''(\tu_0)\tu_1 - \eta_1)(-L_\fil\tilde u_1+\vec{\kappa}\cdot\nabla_\vz U_\fil)+
	\eta_dW''(U_\fil)\tilde u_1.
\eeq
For fixed values of $s$ and $S$, the expression (\ref{MP-eq:t1Tildeu2}) can be solved for $\tilde u_2$ if and only if the right-hand side is perpendicular to
$$ \ker L_\fil = \textrm{span}\{\partial_R U_\fil\cos\theta,\partial_R U_\fil\sin\theta\}=\ker L_{\fil,1}.$$
We decompose $\cQ$ into its~$\cZ_m$ components
\beq
    \cQ = \cQ_0+\cQ_1+\cQ_{0,2},
\eeq
where~$\cQ_0\in\cZ_0$,~$\cQ_1\in\cZ_1$,~$\cQ_{0,2}\in\cZ_0+\cZ_2$, are given by
\begin{align}
    \cQ_0 &:= -W'''(U_\fil)\tilde u_1L_\fil\tilde u_1+\partial_S^2L_\fil\tilde u_1+\eta_1L_\fil\tilde u_1+\eta_dW''(U_\fil)\tilde u_1,\\
    \cQ_1 &:= -\vec{\kappa}\cdot\nabla_\vz L_\fil\tilde u_1+W'''(U_\fil)\tilde u_1\vec{\kappa}\cdot\nabla_\vz U_\fil-\partial_S^2\vec{\kappa}\cdot\nabla_\vz U_\fil-\eta_1\vec{\kappa}\cdot\nabla_\vz U_\fil,\\
    \cQ_{0,2} &:= (\vec{\kappa}\cdot\nabla_\vz)^2U_\fil.
\end{align}
Since the spaces $\cZ_m$ are orthogonal and $L_\fil(\tilde Q_2)\bot\ker L_\fil$,  the solvability conditions take the form
\beq\label{MP-eq:t1SolveCond}
    \left\langle \tilde\mu_2-\cQ_1,\partial_{z_i} U_\fil \right\rangle_{R}=0,\quad \text{ for }i=1,2.
\eeq
To evaluate these conditions we expand~$\cQ_1$ using the expression~(\ref{MP-eq:TildeU1}) for $\tilde u_1$,
\beq\label{MP-eq:t1Q1SolveCond}
	\langle\cQ_1,\partial_{z_i}U_\fil\rangle_{R} = -{2}\pi m_\fil \mu_1\kappa_i-\pi \sigma_\fil\left(\eta_1\kappa_i+ \partial_S^2\kappa_i\right), \quad \textrm{for}\,\, i=1,2,
\eeq
where we have introduced
\begin{align}
    m_\fil&:=\int_0^\infty\hat U_\fil\,RdR,\label{MP-eq:S1}\quad \sigma_\fil:= \int_0^\infty (U_\fil')^2R\,dR.
\end{align}
The  $R$-weight inner product of $\tilde\mu_2$, given in (\ref{MP-eq:t0TildeMu1}), with $\partial_{z_i}U_\fil$ yields
\beq\label{MP-eq:t1mu2SolveCond}
	\langle \tilde \mu_2,\partial_{z_i}U_\fil\rangle_R = \pi V_{i} m_{\fil,2},
\eeq
where we have introduced
\beq
    m_{\fil,2}:=\int_0^\infty\hat U_\fil^2\,RdR\label{MP-eq:S2}.
\eeq
Substituting (\ref{MP-eq:t1Q1SolveCond}) and~(\ref{MP-eq:t1mu2SolveCond}) into (\ref{MP-eq:t1SolveCond}) we arrive at the expression for the normal velocity
\beq\label{MP-eq:bVsS}
	\bV(s,S) = -\frac{2\mu_1m_\fil-\eta_1\sigma_\fil}{m_{\fil,2}}\vkappa-\frac{\sigma_\fil}{m_{\fil,2}}\partial_S^2\vkappa.
\eeq
Introducing the quantities
\beq\label{e:nuf-def2}
\mu_\fil ^*:= \frac{\eta_1\sigma_\fil}{2m_\fil}, \hspace{0.5in} \nu_\fil:= \frac{2m_\fil}{m_{\fil,2}}, \hspace{0.5in} k_\fil = \frac{\sigma_\fil}{m_{\fil,2}},
\eeq
and return the $S$ variable to its original scaling, we obtain the normal velocity
\beq\label{MP-eq:bV}
	\bV(s,S) = -\left[ \nu_\fil(\mu_1-\mu_\fil^*)\vkappa+\eps k_\fil\partial_s^2\vkappa\right].
\eeq

The constant value of $u_1$ is determined by the conservation of total mass, and is coupled to changes in length of the curve $\Gamma_\fil$.
Combining the inner and outer expansions of $u$ yields the composite expansion
\beq\label{MP-eq:t1InnerSol}
    u(x,t) = U_\fil + \eps(\mu_1\Phi_{\fil,2}-\eta_dL_\fil^{-2}W'(U_\fil))+O(\eps^2)\quad\text{in }\Gamma_{\fil,\ell},
\eeq
which has the far-field asymptotics,
\beq\label{MP-eq:t1OuterSol}
    u(x,t) = b_-+\eps\frac{\mu_1}{\alpha_-^2}+O(\eps^2)\quad\text{in }\tilde\Gamma_{\fil,\ell}.
\eeq
The total mass of the system is given by
\beq\label{MP-eq:Mass}
    M:=\int_\Omega u(x,t)-b_-\,dx = \int_\Omega u(x,0)-b_-\,dx= \int_{\Omega \backslash \Gamma_{\fil,\ell}}(u-b_-)\,dx+\int_{\Gamma_{\fil,\ell}}(u-b_-)\,dx,
\eeq
where the outer integral takes the value
\beq\label{MP-eq:t1OuterMass}
	 \int_{\Omega \backslash \Gamma_{\fil,\ell}}(u-b_-)\,dx = \eps \frac{\mu_1}{\alpha_-^2}(|\Omega| -{|\Gamma_{\fil,\ell}|})+O(\eps^2).
\eeq
Using~(\ref{MP-eq:t1InnerSol}) and the Jacobian, (\ref{e:JacobianFil}), we evaluate the inner integral
\begin{align}\label{MP-eq:t1InnerMass}
	\int_{\Gamma_{\fil,\ell}}(u-b_-)\,dx &= \eps^2\int_{\Gamma_\fil}\int_{\mbbR^2}\left(\hat U_\fil+\eps(\mu_1\Phi_{\fil,2}-\eta_dL_\fil^{-2}W'(\phi_\fil))+O(\eps^2)\right)(1-\eps \vz\cdot\vec{\kappa})\,d\vz\,ds\\
	&=\eps^22\pi |\Gamma_\fil |  m_\fil+ O(\eps^3|\Gamma_\fil|).\nonumber
\end{align}
Assuming that $|\Gamma_\fil|\sim O(\eps^{-1})$, as is commensurate with an $O(1)$ amphiphilic mass, we expand
\beq\label{MP-eq:GammaExp}
    |\Gamma_\fil|=\eps^{-1}\gamma_{\fil,-1}+ \gamma_{\fil,0}+ O(\eps),
\eeq
to arrive at the total mass expansion
\beq\label{MP-eq:M}
    M = \eps\left(\frac{\mu_1}{\alpha_-^2}|\Omega|+2\pi m_\fil\gamma_{\fil,-1}\right)+O(\eps^2).
\eeq
Taking the~$\tau=\eps t$ time derivative of the total mass,~(\ref{MP-eq:M}), and solving for~$\frac{d\gamma_{\fil,-1}}{d\tau}$ yields
\beq\label{MP-eq:ChangeMu1}
	\frac{d\gamma_{\fil,-1}}{d\tau} = -\frac{|\Omega|}{2\pi\alpha_-^2m_\fil}\frac{d\mu_1}{d\tau}.
\eeq
On the other hand, any admissible codimension two curve evolving with normal velocity $\bV$ satisfies
\beq\label{MP-eq:t1GammaGrowth}
    \frac{d|\Gamma_\fil|}{d\tau}= -\int_{\Gamma_\fil}\bV\cdot\vec\kappa\,ds.
\eeq
Combining this expression with (\ref{MP-eq:bV}), (\ref{MP-eq:GammaExp}), and (\ref{MP-eq:ChangeMu1}) yields
\beq\label{MP-t1mu1ODE}
	\frac{d\mu_1}{d\tau} = \eps \frac{2\pi \alpha_-^2 m_\fil}{|\Omega|} \left(-\nu_\fil \left(\mu_1-\mu_\fil^*\right)\int_{\Gamma_\fil}|\vec\kappa|^2\,ds+
     \eps k_\fil\int_{\Gamma_\fil}|\partial_s\vec\kappa|^2\,ds\right) + O(\eps^2).
\eeq
This system exhibits the same quenching behavior as the codimension one evolution, with the distinction being
that the equilibrium far-field density for a codimension two curve takes the form
\beq
    \lim\limits_{\tau\to\infty} u = b_- +\eps\frac{\mu_\fil^*}{\alpha_-^2}+O(\eps^2).
\eeq

\section{Competitive evolution of amphiphilic suspensions}

We assume that $\Omega\subset\RR^3$ and combine the geometric flow results derived in sections \ref{Sec:Bilayer} and \ref{Sec:Filament} with the pearling stability results for
bilayers and filaments presented in \cite{NK-KP-18}.  The goal is to derive an overall picture of the complexity of transients and bifurcation structure of the $H^{-1}$
gradient flow of the strong scaling of the FCH system.

\subsection{Competitive evolution of codimension one and two systems}\label{Sec:BilayerVsPore}

Fix $\Omega\subset\RR^3$ and let $\Gamma_b$ and $\Gamma_\fil$ be admissible codimension one and codimension two morphologies with disjoint reaches, $\Gamma_{b,\ell}$ and $\Gamma_{\fil,\ell}$. We emphasize that $\Gamma_b$ and $\Gamma_\fil$ may be comprised of multiple disjoint surfaces and curves.
For a given value of the chemical potential, $\mu_1$, the codimension one bilayer morphology $u_b$, given in (\ref{MB-eq:t1InnerSol}) and
codimension two filament morphology $u_f$, given in (\ref{MP-eq:t1InnerSol}) satisfy identical far-field asymptotics
\beq
	\lim\limits_{R\to\infty} u_{\fil}=\lim\limits_{z\to\infty}u_{b} = b_- -\eps\frac{\mu_1}{\alpha_-^2}+O(\eps^2).
\eeq
Consequently we may form the composite solution
\beq\label{e:Composite}
    u_{b,\fil} = u_b+u_\fil -\left(b_- -\eps\frac{\mu_1}{\alpha_-^2}\right)+O(\eps^2),
\eeq
parameterized by $\Gamma_b$, $\Gamma_\fil$, and the common, slowly varying, chemical potential~$\mu_1$.
Recalling the scalings (\ref{MB-eq:GammaExp}) and (\ref{MP-eq:GammaExp}) of the surface area and length of $\Gamma_b$ and $\Gamma_\fil$ respectively,
the total mass of the composite solution satisfies
\beq
	M =  \eps\left(\frac{\mu_1}{\alpha_-^2}|\Omega| +  m_b\gamma_{b,0}+ 2\pi m_\fil \gamma_{\fil,-1}\right)+O(\eps^2),
\eeq
where $m_b$, the bilayer mass per unit area, is defined in (\ref{e:MB-const}), and $ 2\pi m_\fil$ denotes the filament mass per unit length, defined in
(\ref{MP-eq:S1}). Expanding $M=\eps \hat M+O(\eps^2)$ and solving for $\mu_1$ yields relation between the morphology size and the chemical potential~$\mu_1$,
\beq\label{MB-eq:Mu1Constraint}
	\mu_1 = \frac{\alpha_-^2}{|\Omega|}\left(\hat M-m_b\gamma_{b,0} -  2\pi m_\fil\gamma_{\fil,-1} \right).
\eeq
Taking the time derivative of~(\ref{MB-eq:Mu1Constraint}) and using the relations (\ref{MB-eq:gammaGrowthBl}) and
(\ref{MP-eq:t1GammaGrowth}) to relate
the growth of the curves to the normal velocities, yields an evolution equation for the chemical potential
given in (\ref{e:Mu1Evol}).
%
Coupling this equation to the normal velocities for the bilayer and filament derived in (\ref{e:NVb}) and (\ref{e:NVf}) gives a closed system for
the combined curve motion and far-field chemical potential.

\begin{figure}[h!]
\begin{center}
\begin{tabular}{ccp{2.5in}}
    \includegraphics[width=2.05in]{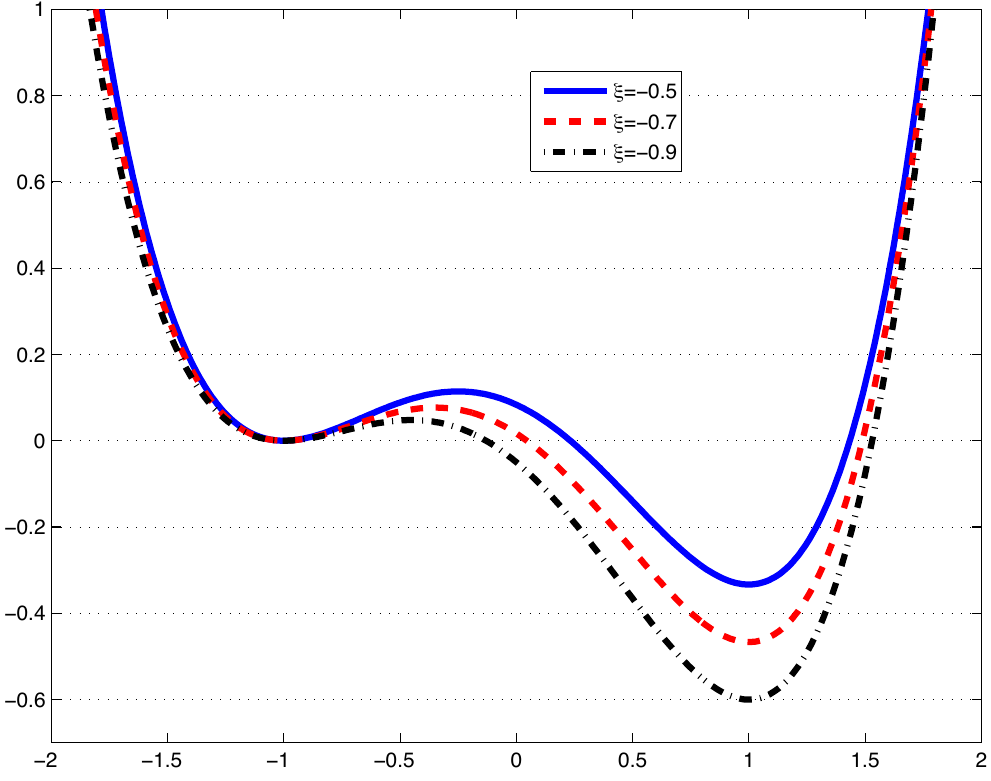} &
    \includegraphics[width=2.05in]{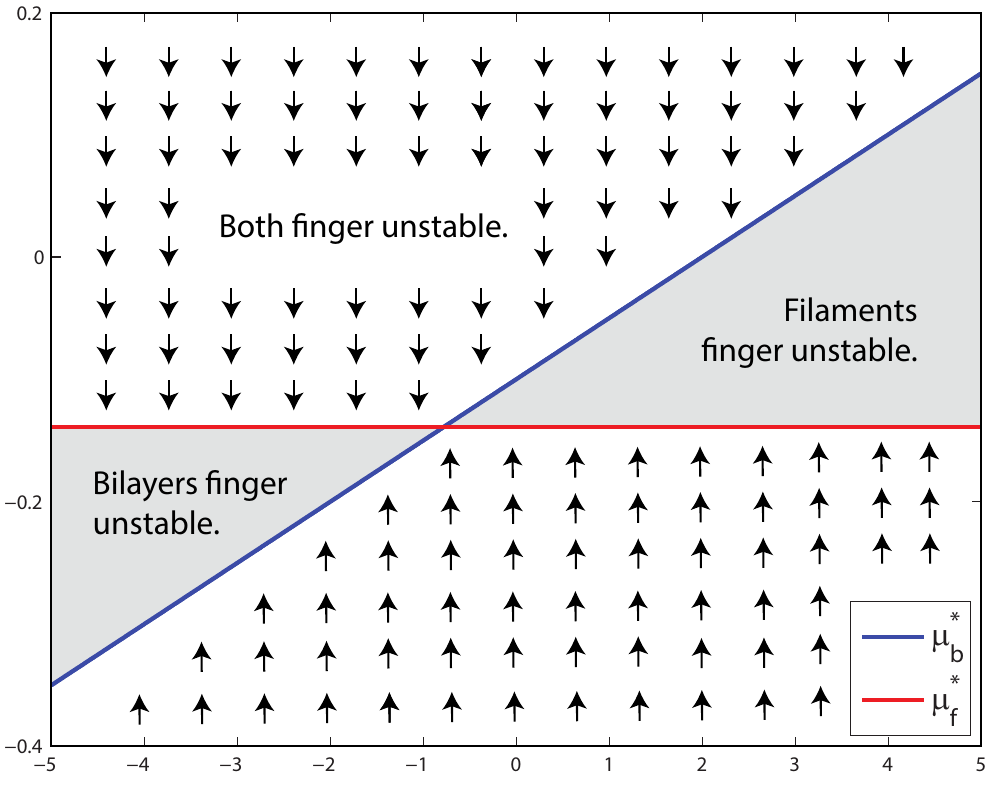}&
    \includegraphics[width=2.06in]{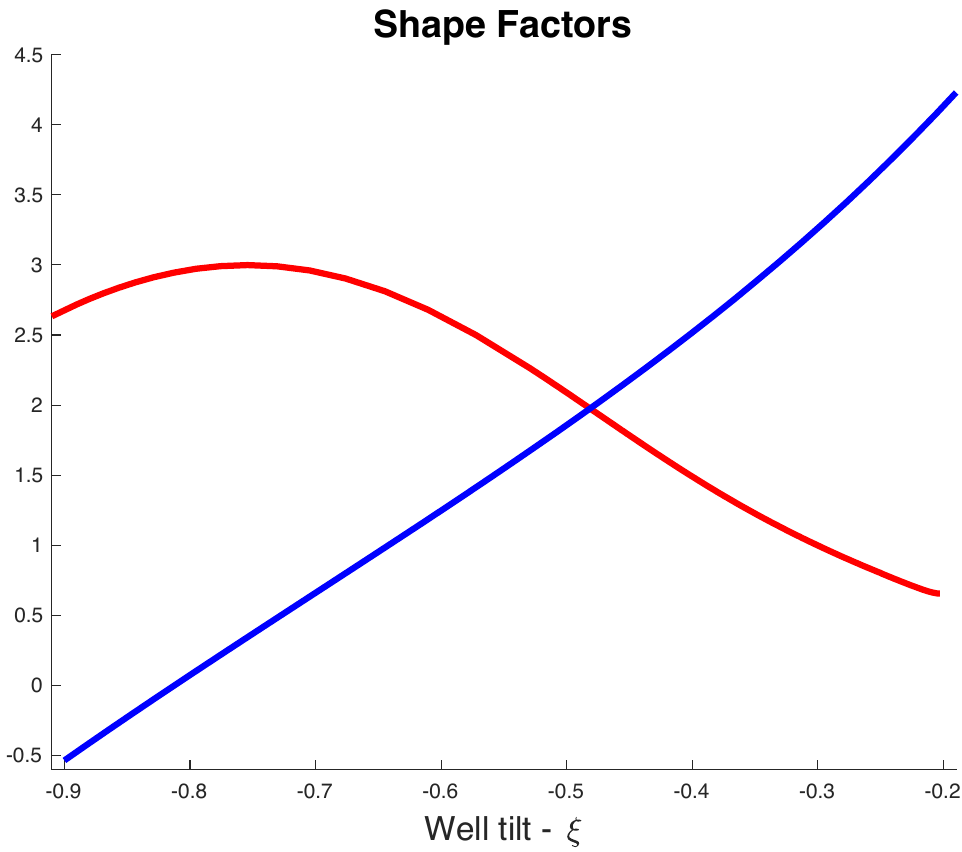}
\end{tabular}
 \vskip 0.1in
  \caption{(left) Graph of the double well, $W$ versus $u$ for $\xi=-0.9,-0.7,-0.5$ (see (\ref{NV-eq:W})). (center) Diagram of the curve shortening regions in $\eta_d-\mu_1$ plane, for
 $\xi=-0.9$ and $\eta_1=0.15$. Arrows indicate the direction of the temporal evolution of $\mu_1$ under the flow (\ref{e:Mu1Evol}).
 The shaded region is forward invariant and globally attracting so long as the curves remain admissible. (right) Values of the bilayer and filament shape factors, $S_b$ (blue) and $S_\fil$ (red) as a function of the well tilt, $\xi$. The change in sign of $S_b$ near $\xi=-0.8$ flips the direction of the inequality in (\ref{e:BPS}).}
 \label{f:W-xi}
 \end{center}
\end{figure}
\subsection{Analysis of competitive geometric evolution and bifurcation}

The analysis presents a bifurcation diagram with four thresholds that delineate distinct behaviors. The $\mu_1$ thresholds for the pearling bifurcation, $P_b$ and $P_\fil$, depend upon the functionalization parameters  $\eta_1$ and $\eta_2$ through their difference $\eta_d:=\eta_1-\eta_2.$  The $\mu_1$ thresholds for the transition from curve shortening to regularized curve lengthening, $\mu_b^*$ and $\mu_\fil^*$ given in (\ref{e:nub-def2}) and (\ref{e:nuf-def2}}) respectively, have a more subtle dependence upon the functionalization parameters. These
relations are summarized below:
\begin{align}
\label{e:BPS}
\mu_1\,\textrm{sign}(S_b)<  P_b(\eta_d) &:= -\eta_d \frac{\lambda_{b,0} \|\psi_{b,0}\|_{L^2}^2}{|S_b|}, && \textrm{Bilayers Pearling Stable} \\
\label{e:FPS}
\mu_1\,\textrm{sign}(S_\fil) < P_\fil(\eta_d)&:= - \eta_d \frac{\|\psi_{\fil,0,0}^\prime\|_{L^2_R}^2 +\lambda_{\fil,0,0}\|\psi_{\fil,0,0}\|^2_{L^2_R}}{|S_f|}, && \textrm{Filaments Pearling Stable} \\
\label{e:BFS}
&\mu_1<\mu_b^*(\eta_1,\eta_2), && \textrm{Bilayers Curve Shortening}\\
\label{e:FFS}
&\mu_1<\mu_\fil^*(\eta_1,\eta_2),&&\textrm{Filament Curve Shortening}.
\end{align}
The signs of the shape factors $S_b$ and $S_\fil$, defined in (\ref{eq:SF-def}), depend upon the choice of the double well, $W$, and impact not only the sign
 of the right-hand sides of (\ref{e:BPS}) and (\ref{e:FPS}) but also the direction of the inequalities, see Figure\,\ref{f:W-xi} (right). The chemical potential and the geometric flows evolve on the same $t=O(\eps^{-1})$ timescale. Within the $H^{-1}$ gradient flow, the pearling instability produces eigenvalues of
size $O(\eps^{-1})$, \cite{NK-KP-18} and hence will manifest itself on the $t=O(\eps)$ timescale, essentially instantaneously on the time scale of the geometric flow and the chemical potential. We define the \emph{pearling instability region} to be set of values $(\mu_1,\eta_1,\eta_2)$ for which either codimension one
or codimension two structures are pearling unstable.


We investigate the variation of the pearling instability regions, and its relation to regions of curve lengthening flows, as functions of the well shape.
For simplicity we present these regions in the $\mu_1-\eta_d$ plane, with the assumption that $\eta_1=0.15$, unless specified otherwise. To parameterize the well
shape we fix $b_\pm=\pm1$ and insert a one-parameter well tilt, $\xi$ into the double-well potential,
\beq\label{NV-eq:W}
   W(u;\xi) := \frac{1}{2}(u-b_-)^2\left(\frac12(u-b_+)^2-\frac{\xi}{3}\left(u-\frac{3b_+-b_-}{2}\right)\right),
\eeq
where the parameter~$\xi$ controls the value of $W$ at the right well $u=b_+$, see Figure\,\ref{f:W-xi} (left).

\begin{figure}[h!]
\begin{center}
  \begin{tabular}{cccc}
        \includegraphics[width=1.5in]{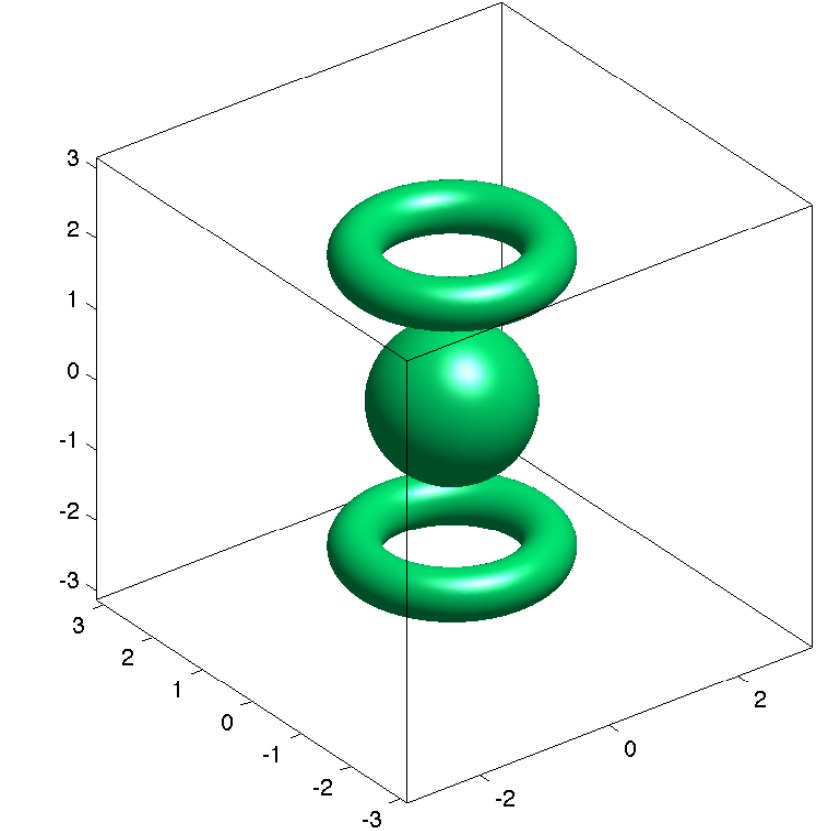} &
      \includegraphics[width=1.5in]{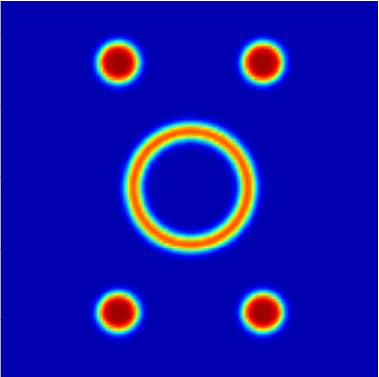} &
         \includegraphics[width=1.5in]{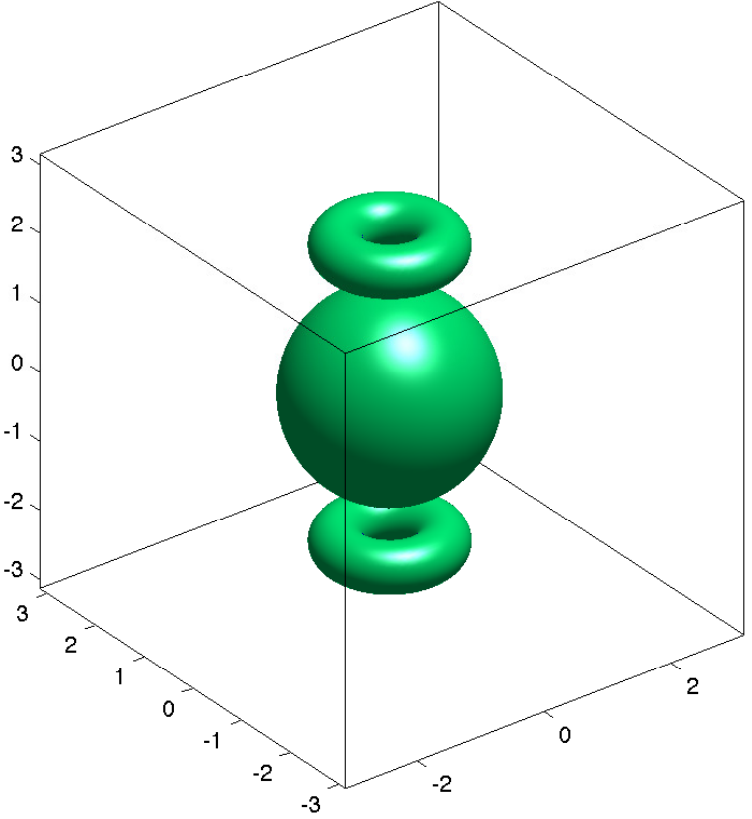} &
      \includegraphics[width=1.5in]{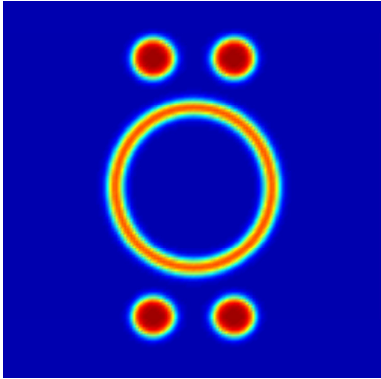} \\
      (a) $t=0$ & (b) $t=0$ slice & (c) $t=200$ & (d) $t=200$ slice
   \end{tabular}
  \end{center}
  \vskip -0.2in
\caption{Simulation of the FCH gradient flow for $\eps=0.03$, $\eta_1=0.15$, $\eta_2=0.24$, and $\xi=-0.15$. For these values $\mu_b^*<0<\mu_\fil^*$ and after a short
transient bilayers will grow while filaments
shrink. (a) Initial data consisting of two circular filaments and a spherical bilayer inside the $[-\pi, \pi]^3$ computational domain. (b) 2D slice along $x-z$ plane ($y=0$) of initial data.
(c)  3D pose of the $t=200$ $(\tau=6)$ final computation stage and (d) corresponding 2D slice showing the larger bilayer and smaller filament radii.
}
 \label{f:T200}
\end{figure}

For the $\mu_1$ dynamics, as $\eta_1$ and $\eta_2$ are fixed parameters, the temporal evolution of $\mu_1$ traces a vertical line
segment on the diagram, with $\mu_1$ decreasing if it is larger than both $\mu_b^*$ and $\mu_\fil^*$ and increasing if it is smaller than both.
The $\mu_1$-region bounded by the points $\{\mu_b^*,\mu_\fil^*\}$ is attracting and forward invariant under the flow
see Figure\,\ref{f:W-xi} (center). Within this invariant region the direction of the flow depends upon the overall
curvatures of the two classes of morphologies; however to leading order the total area/length of the codimensional structures with the larger equilibrium
value decreases monotonically and the other increases monotonically, so long as the morphologies remain admissible.
A simulation of the FCH gradient flow with one spherical bilayer and two circular filaments is presented in Figure\,\ref{f:T200} for parameter values for
which $\mu_b^*<\mu_\fil^*$.

\begin{figure}[h!]
\begin{center}
  \begin{tabular}{ccc}
        \includegraphics[width=2.1in]{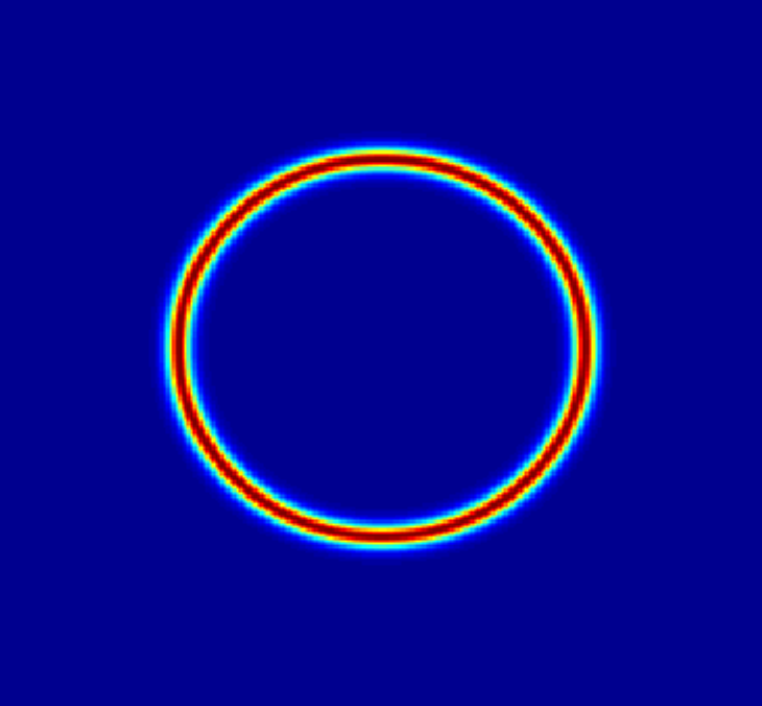} &
      \includegraphics[width=2.1in]{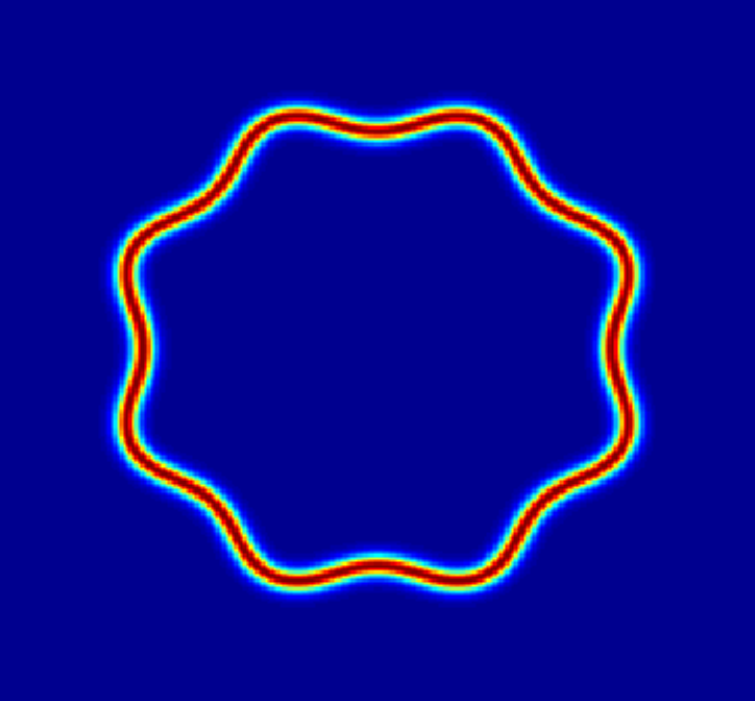} &
      \includegraphics[width=2.1in]{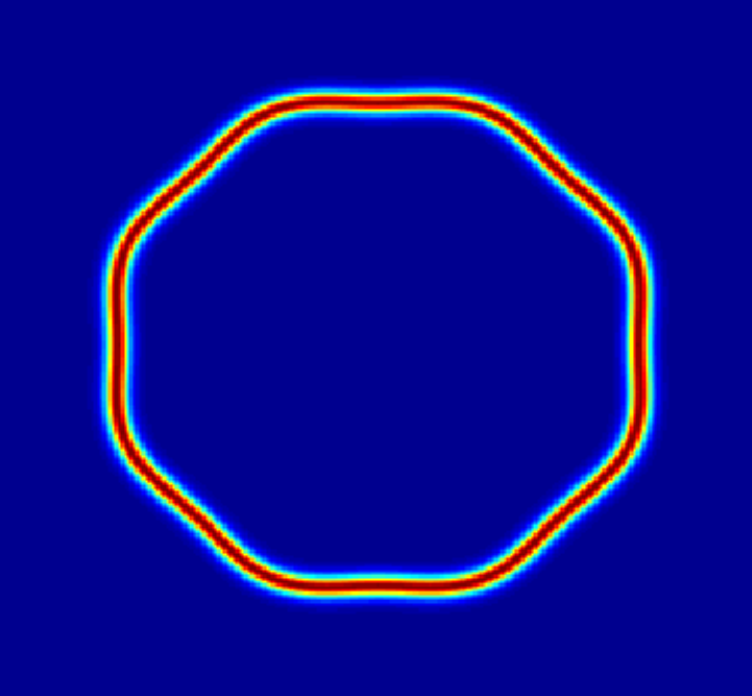} \\
     (a) $t=0$ & (b) $t=484$ & (c) $t=986$
   \end{tabular}
  \end{center}
  \vskip -0.3in
\caption{Simulation of the FCH gradient flow corresponding to a single circular bilayer in $[-2\pi, 2\pi]^2\subset \RR^2$  in the regularized curve lengthening regime. The larger
radii and value of $\mu_1-\mu_b^*$, compared to simulations presented in Figure\,\ref{f:T200}, induce the interfacial meander. Yet larger initial values of $\mu_1-\mu_b^*$ lead to
self-intersection and defect generation. (a) The initial data has $\mu_1>\mu_b^*$, with sufficient excess to initiate shape perturbations as the circular bilayer grows. (b)  The onset of the shape instability at time $t=484$. (c) Quenching of the flow at $t=986$ as the value of $\mu_1$ relaxes towards $\mu_b*$ and the higher order surface diffusion returns the interface back to a larger, circular shape.
}
 \label{f:Transman}
\end{figure}

Motivated by this example, and to further illustrate the dynamics, we
consider a composite morphology consisting of $N_b$ spherical bilayers of radii $R_1, \ldots, R_{N_b}$ and $N_\fil =O(\eps^{-1})$
circular filaments of radii $r_1, \ldots, r_{N_\fil}$. \tck{The set of such morphologies is forward invariant under the gradient flow assuming that the
pearling stability conditions hold. However if the flow resides in either of the regularized curve lengthening regimes ($\mu_1> \mu_b^*$ or $\mu_1>\mu_\fil^*)$ then
the spherical or circular shapes, respectively, are not stable under this flow.} For these special shapes the competitive evolution (\ref{e:NVb})-(\ref{e:Mu1Evol}) reduces to
\begin{align}
\dot{R}_i &= \nu_b\left(\mu_1 -\mu_b^*\right) \frac{2}{R_i}, \hspace{0.5in} i=1,\ldots, N_b, \label{e:R-evol}\\
\dot{r}_j&=\nu_\fil \left(\mu_1-\mu_\fil^*\right) \frac{1}{r_j },\hspace{0.5in} j=1,\ldots, \eps^{-1}N_\fil, \label{e:r-evol}\\
\dot{\mu}_1 &= -\frac{\alpha_-^2}{|\Omega|} \left( 16\pi  m_b\nu_b (\mu_1-\mu_b^*)N_b +\eps4\pi^2 m_\fil \nu_f(\mu_1-\mu_f^*)\sum\limits_{j=1}^{N_\fil} \frac{1}{r_j}\right),\label{e:mu-evol}
\end{align}
where the dot notation denotes differentiation with respect to $\tau=t/\eps.$
The $\mu_1$ evolution depends upon the spherical bilayers only through their total number. Consistent with the discussion above,
 the bulk chemical potential $\mu_1$ decreases if $\mu_1>\max\{\mu_{b}^*,\mu_{\fil}^*\}$ and increases if $\mu_1<\min\{\mu_{b}^*,\mu_{\fil}^*\}$. The radii shrink
 or grow depending upon the signs of $\mu_1-\mu_b^*$ and $\mu_1-\mu_\fil^*$.  If both $N_b$ and $N_\fil$ are positive, and $\mu_b^*\neq\mu_\fil^*$ then the system has no equilibrium.
  Assuming for simplicity of presentation that $\mu_b^*<\mu_\fil^*$, then after a possible transient the system
 enters a regime in which all spheres are growing and all circular filaments are shrinking. The radii $r_j$ will then decrease to zero in finite time due to the
 inverse relation between $\dot{r}_j$ and $r_j$.
 If the zero radius filaments are removed from the system and the remaining filaments re-indexed, then after a transient the set of hoops will be empty ($N_\fil=0$) and only spheres will
 remain. At this point $\mu_1$ will relax (quench) at an exponential rate from above to $\mu_b^*$ as the spherical radii grow, albeit at an exponentially decreasing rate as $\mu_1$
 quenches to $\mu_b^*$. The equilibrium will be a collection of spherical bilayers of differing radii. We emphasize that in the growing regime spherical shapes are unstable
 to perturbation under the full flow (\ref{e:NVb}). Manifestation of this instability requires sufficiently large values of $\nu_b(\mu_1-\mu_b^*)>0$ in relation to the coefficient
 $\eps k_b$ of the surface diffusion term in (\ref{e:NVb}), and is triggered more easily with increasing radius. However, if the non-spherical excursions are not so large as to induce self-intersection,
 then the shapes return to spherical as $\mu_1$ quenches to $\mu_b^*$.  Figure\,\ref{f:Transman} presents a  simulation of the full FCH gradient flow that illustrates this phenomena,
 while a rigorous derivation of the transient instability regime of the curve lengthening flow in $1+2D$ for nearly circular bilayers is presented in \cite{CP-19}.

\begin{figure}[h!]
\begin{center}
  \begin{tabular}{cc}
        \includegraphics[width=3.1in]{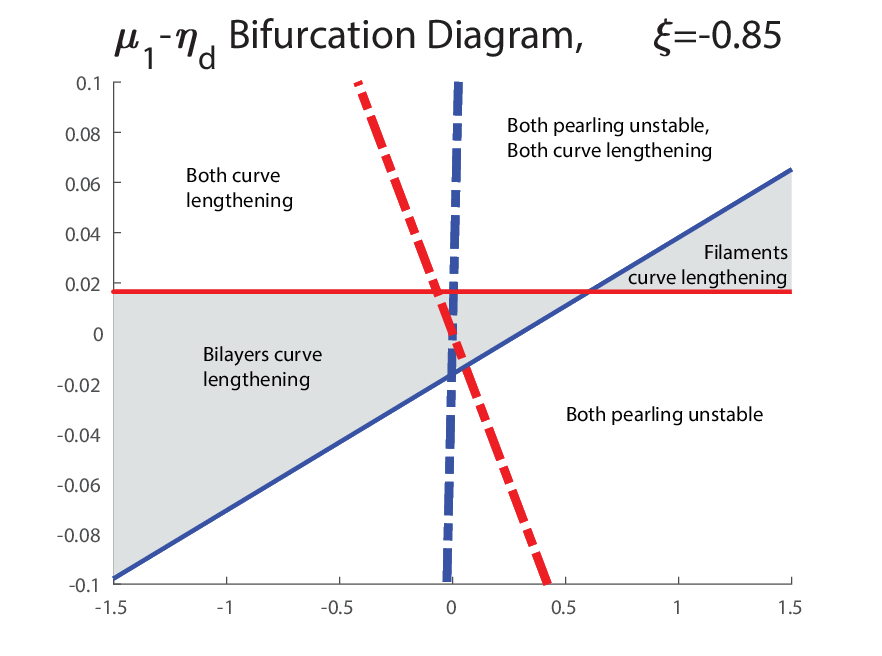} &
      \includegraphics[width=3.1in]{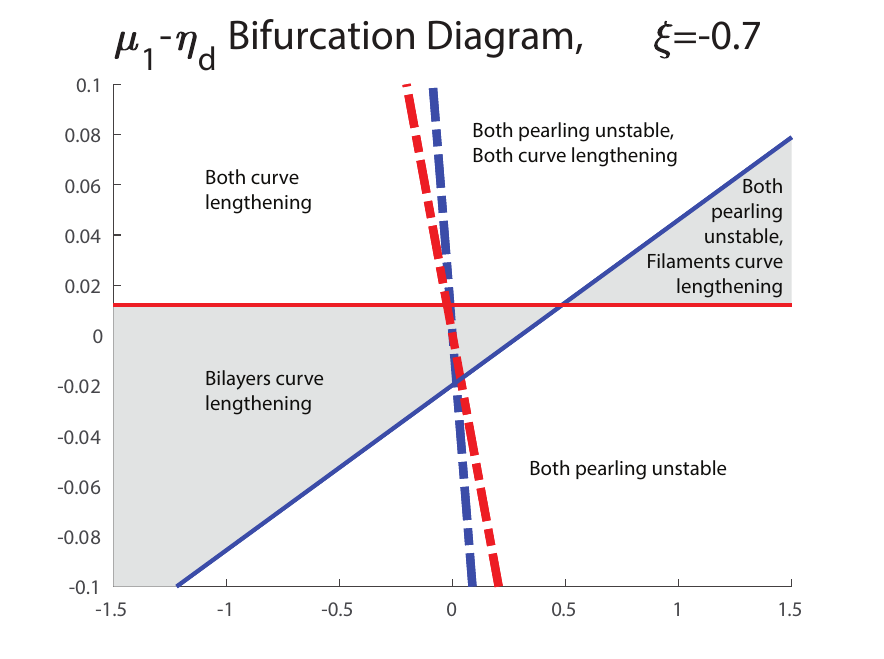} \\
      \includegraphics[width=3.1in]{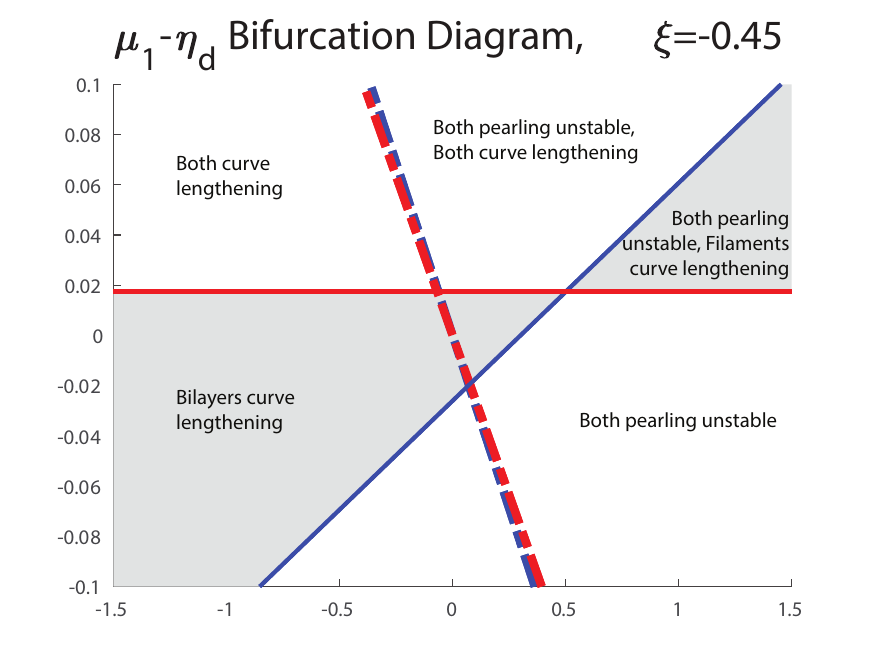} &
       \includegraphics[width=3.1in]{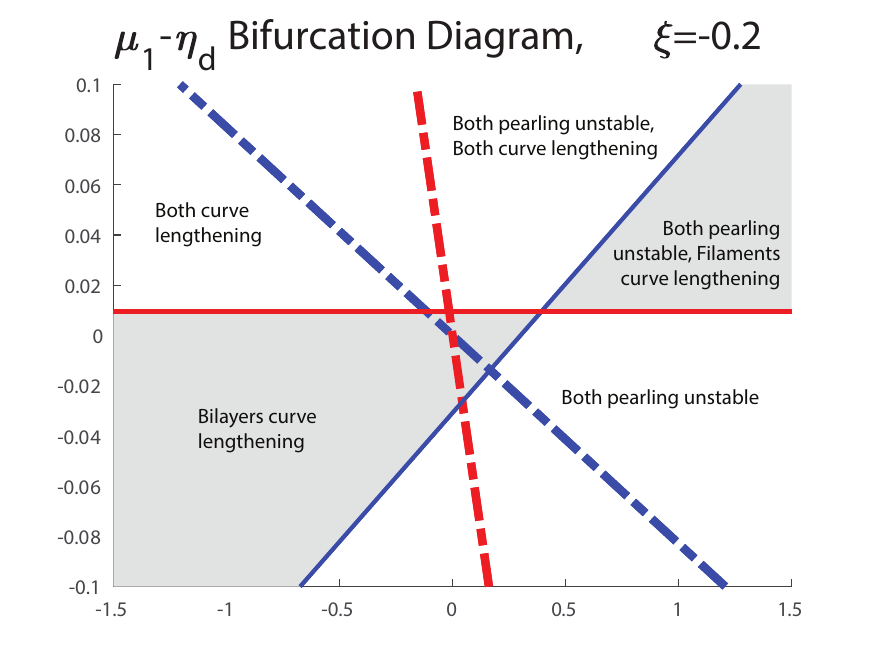}
   \end{tabular}
  \end{center}
\caption{A $\mu_1\!-\!\,\eta_d$ bifurcation diagram versus $\mu_1$ (vertical axis) and $\eta_d$ (horizontal axis) for $\xi=-0.85, -0.7, -0.45$ and $-0.2$ and $\eta_1$ fixed at $\eta_1=0.15$. Pearling instability holds to the right of the dotted lines (blue-bilayers, red-filaments) while regularized curve lengthening holds above the solid lines (blue-bilayers, red-filaments). Areas of pearling instability and curve lengthening are indicated in key regions. The value of $\mu_1$ is generically time dependent, and the grey shaded region is forward invariant
under the flow.
}
 \label{f:Bdiag}
\end{figure}

In Figure\,\ref{f:Bdiag}, the pearling bifurcation lines are added and the full stability diagram is shown for values of the well-tilt parameter
$\xi=-0.85, -0.7, -0.45,$ and $-0.2$.
For $ \xi=-0.85$, $S_\fil$ is positive, $S_b$ is slightly negative, and for positive values of $\mu_1$ filament pearling instability region contains the bilayer pearling
instability region. As $\xi$ is decreased from $-0.7$ to $-0.45$ the two pearling instability lines almost coincide, and for $\xi=-0.2$ they have crossed, with the
filament pearling instability region now lying above the bilayer pearling instability region for $\mu_1>0$.   In all figures the intersections of the pearling instability
lines occurs to the left of the crossing of the curve lengthening lines. This implies that increasing $\eta_d$ from negative values will excite
pearling bifurcations at smaller values of $\eta_d$ than those for which the filament morphology becomes dynamically favored over the bilayer morphology.
For the non-generic value of $\eta_d$  at which the critical values $\mu_b^*=\mu_\fil^*$ coincide, codimension one and two structures can co-exist on the $t=O(\eps^{-1})$
time-scales under consideration here.

\subsection{Analytic bifurcation diagrams and comparison to simulations}
\label{s:ABD-sim}
\begin{figure}[h!]
         \centering
           \begin{subfigure}[b]{0.35\textwidth}
         \includegraphics[width=2.5cm]{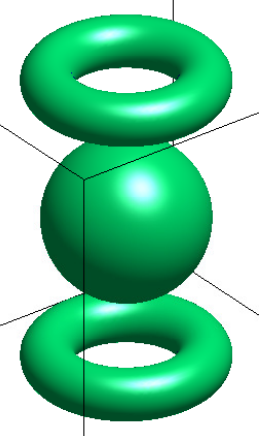} \hspace{0.1in}
        \includegraphics[width=2.5cm]{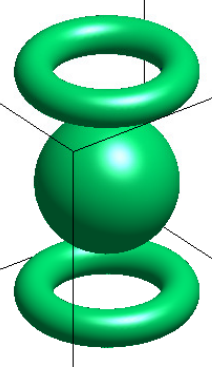} \\
        \includegraphics[width=2.5cm]{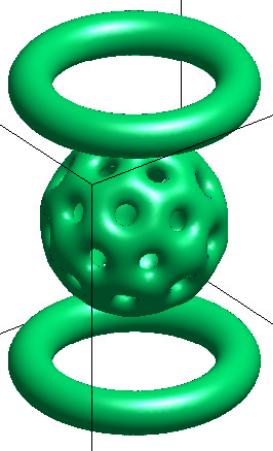} \hspace{0.1in}
        \includegraphics[width=2.5cm]{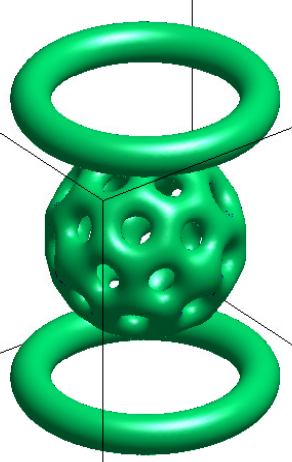}
        \end{subfigure}
        \begin{subfigure}[b]{0.5\textwidth}
           \includegraphics[width=2.9in]{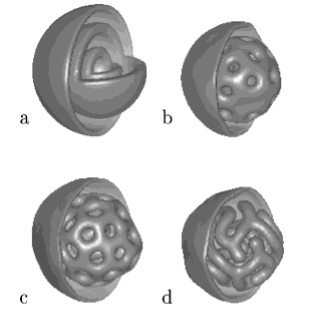}
        \end{subfigure}
       \begin{subfigure}[b]{0.10\textwidth}
         \hspace{1.0in}
          \end{subfigure}
        \vskip -2.8in
        \begin{tabular}{p{4.8in}c}
        &     \includegraphics[width=1.2in]{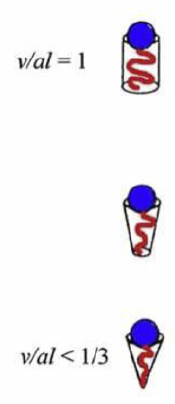}
        \end{tabular}
                \caption{(left) The $t=100$ ($\tau=3$) end state of simulations of the FCH gradient flow
for $\eps=0.03$, $\eta_1=0.15$, $\eta_2=0.24$ from initial data described in the text and well tilt $\xi= -0.15, -0.20$ (top row, left to right),
$\xi= -0.25, -0.3$ (bottom row, left to right). The less negative values of $\xi$ have initial values of $\mu_1$ that inhibit pearling,
for the more negative two values the bilayer pearls. Images courtesy of Andrew Christlieb and Jaylan Jones.
(center) Simulations of a mean-field density functional model of amphiphillic diblock copolymers with ratio of amphiphilic component of the diblock decreasing from (a) 35\% (b) 30\% (c) 25\% (d) 20\%. The minority solvent-hydrophilic phase is imaged.
Reprinted (adapted) with permission from \cite{fraaije2003model}. Copyright 2003 American Chemical Society.
(right) Depiction of aspect ratio of lipids, a biological diblock with a short amphiphilic head and a long hydrophobic tail. The aspect ratio is defined as the lipid volume $v$ divided by the head cross-sectional area $a$ and tail length $l$. From  \cite{mouritsen2011lipids}, reprinted with permission from the Wiley publication. Copyright \textcopyright 2011 WILEY-VCH Verlag GmbH \& Co. KGaA, Weinheim.}
 \label{f:Jaylan}

\end{figure}

We compare the evolution of numerical simulations of the $H^{-1}$ gradient flow of the FCH free energy, (\ref{e:FCH-eq}),  to the corresponding bifurcation results
and to simulations from a self-consistent mean field model. The parameters impacting stability are the functionalization parameters
$\eta_1$ and $\eta_2$, the time-dependent scaled chemical potential $\mu_1$, and the  shape of the double well, parameterized by $\xi.$  For simplicity we fix $\eta_1=0.15$ and
vary $\eta_2$ and $\xi$.  Simulations of the strong FCH gradient flow, (\ref{e:FCHFE}), were conducted in a domain $\Omega=[-\pi,\pi]^3\subset \RR^3$, with initial data of the
form (\ref{e:Composite}) with $\Gamma_b$ consisting of a single sphere of radius $R_1=1.1$ with center
at $(0,0,0)$ and $\Gamma_\fil$ comprised of two circular filaments of radius $r_1=r_2=1.6$ oriented parallel to the $x-y$ plain of $\Omega$
and with center point at $\left(0,0,\pm\frac12(\pi+R_1/2)\right).$ The initial value of $\mu_1$ varied with each simulation  and is reported in
Figure\,\ref{f:Jaylan}. The system parameters are  $\eps =0.03$, $\eta_1=0.15$ and $\eta_2=0.24$, hence $\eta_d=-0.09$, and four separate values of
$\xi:=-0.15, -0.2, -0.25, -0.3.$ In addition, there was one simulation for $\xi=-0.2$ and $\eta_1=\eta_2=0.15.$
The simulations were conducted for $t\in[0,100]$, equivalently $\tau\in[0,3]$.

\begin{figure}[h!]
\begin{center}
  \begin{tabular}{cc}
      \includegraphics[width=3.1in]{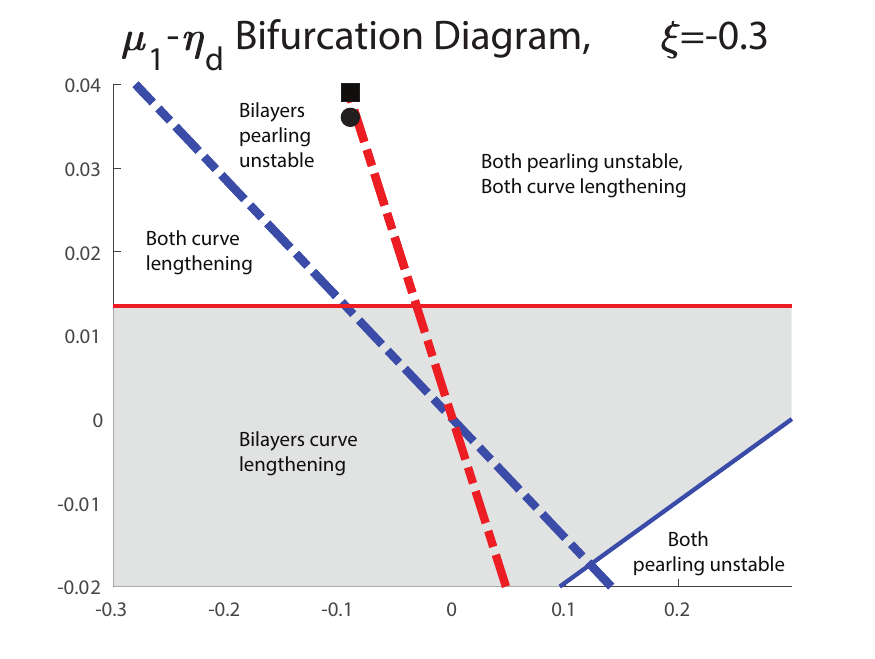} &
      \includegraphics[width=3.1in]{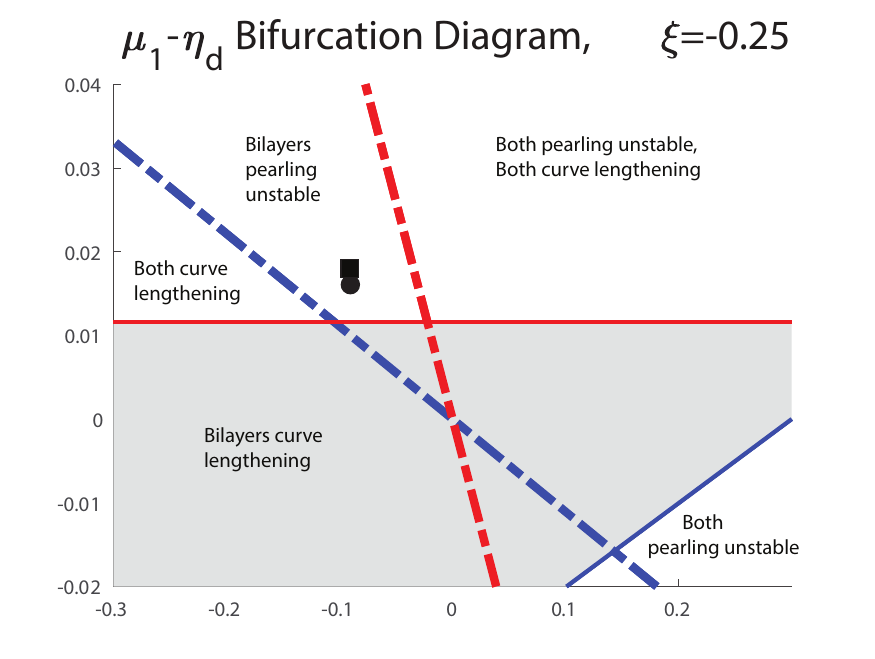} \\
      \includegraphics[width=3.1in]{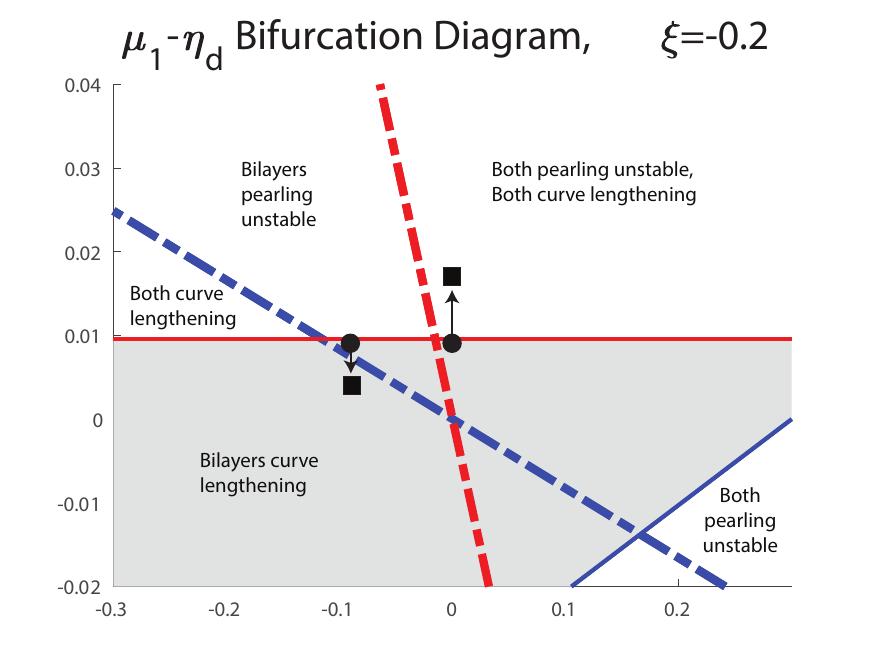} &
      \includegraphics[width=3.1in]{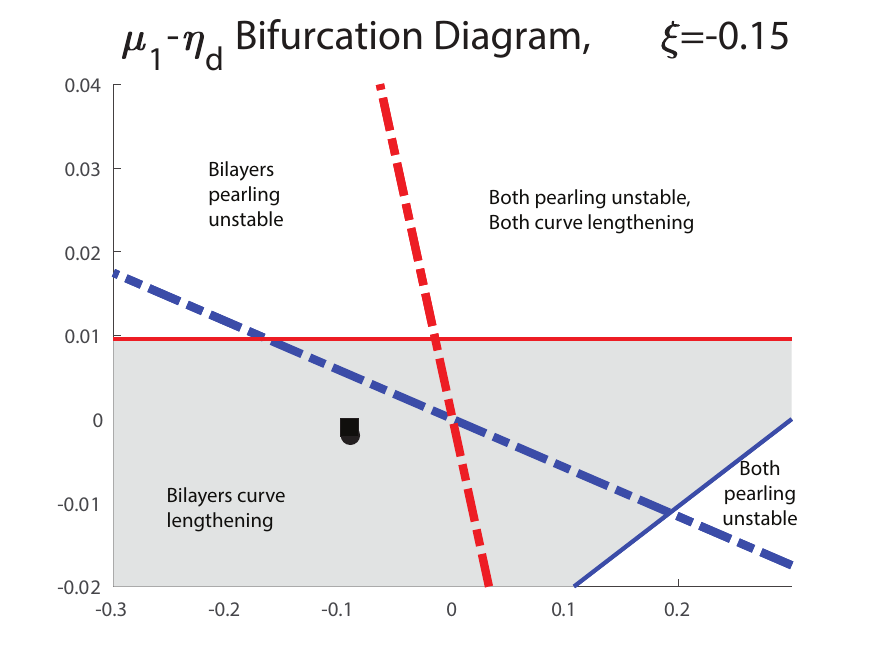} \\
  \end{tabular}
  \end{center}
\caption{$\mu_1\!-\!\eta_d$ bifurcation diagrams with $\eta_1=0.15$ for the same values of $\xi$ as in the end-states depicted in Figure\,\ref{f:Jaylan} (left). For the five simulations the initial and $t=100$ final value of $\mu_1$ is indicated on the corresponding bifurcation diagram
with a solid black circle and square, respectively. The simulations with $\xi=-0.25, -0.3$ and $\xi=-0.2$ with $\eta_d=0$ lead to pearled morphologies.}
 \label{f:BdiagJ}
\end{figure}

Results for the five simulations are superimposed upon the corresponding $\mu_1\!-\!\eta_d$ bifurcation diagram and presented in Figure\,\ref{f:BdiagJ}.
The initial and final values of $\mu_1$ for each simulation are indicated with a closed circle and closed square respectively in each of the four
diagrams. For $\xi=-0.15$ the value of $\mu_1$ starts in a region of bilayer and filament pearling stability, filament curve shortening,
and bilayer regularized curve lengthening, see Figure\,\ref{f:BdiagJ}(bottom-right). During the simulation the bilayer radius grew, the filaments shrunk, and
neither pearled, the $\tau=3$ end state is presented in Figure\,\ref{f:Jaylan} (left/top-left). Two simulations were conducted for $\xi=-0.2$,
for the simulation with $\eta_d=-0.09$ the initial value of $\mu_1$ lies at the boarder of the bilayer pearling region, and the initial stages
of the simulation $(t<5)$ displayed the onset of pearling, however the value of $\mu_1$ decreased out of the pearling region as the filaments
and bilayers grow and the pearling evanesced, restoring the unpearled bilayer structure. The end-state is presented in
Figure\,\ref{f:Jaylan} (left/top-right), there is less shrinking of the circular filaments than in the case $\eta_2=-0.15$ and
the filaments are thinner due to the stronger well tilt. The $\xi=-0.2$ simulation with $\eta_d=0$ but an identical initial value
of $\mu_1$ starts in the middle of the bilayer pearling region, the bilayer pearled fully (end-state not shown) and the value of $\mu_1$ increased,
see Figure\,\ref{f:BdiagJ} (bottom-left). The simulations with $\xi=-0.25$ and $\xi=-0.3$ begin
within the bilayer pearling region, see Figure\,\ref{f:BdiagJ} (top-left and top-right) and rapidly pearled with the $\xi=-0.25$ simulation
pearling around $t=5$ and the $\xi=-0.3$ simulation fully pearled at the first output time of $t=1$.  The pearling  leads to an
increase in $\mu_1$ as the bilayer sheds net amphiphilic molecule mass to the bulk (far-field). The end states are depicted in
Figure\,\ref{f:Jaylan} (left/bottom-left and bottom-right).

 The bilayer pearling instability was observed by Fraaije and Sevink, who developed a self-consistent mean field density
functional model describing the free energy of amphiphilic diblock polymer surfactants embedded in solvent. Their model parameters are based
upon  poly(propylene oxide)-poly(ethylene oxide) diblock in an aqueous solution, see \cite{fraaije2003model} and reference therein.
They simulated spherical nanodroplets of 15\% solvent and 85\% polymer by volume. By decreasing the block ratio -- the ratio of
the length of the hydrophilic portion of the diblock chain to the length of the hydrophobic portion, they uncovered a series of bifurcations that
lead from stable bilayers, to pearled bilayers to a continuous filament pattern decorated with $Y$-junctions and endcap defects, see Figure
\ref{f:Jaylan} (center). This change in block ratio from $35$\% amphiphilic polymer down to $20$\% amphiphilic polymer increases the
aspect ratio of the minority phase.
Within the context of the FCH,  the parameter $\eta_2$ weighs the energy requirement of compressing
the minority phase into the restricted core of a higher codimensional structure, corresponding to values of $u$ for which the double well $W$ is negative.
In particular, bilayer profiles reside in the positive region of the double well, $W(\phi_b)>0$, while the filament profile accesses the negative regions of
$W$, consequently negative values of $\eta_2$ lower the energy cost of filaments relative to bilayers. In this sense positive values of
$\eta_2$ penalize higher codimensional structures, which is analogous to amphiphilic molecule aspect ratio near one, while negative values of $\eta_2$ favor the packing
found in high codimension morphologies, analogous to a high aspect ratio amphiphilic molecule, see Figure\,\ref{f:Jaylan} (right). These preferences
are born out in Figure\,\ref{f:Bdiag}: for fixed values of $\mu_1$, decreasing $\eta_2$ under constant $\eta_1$, hence increasing $\eta_d$, results in
a crossing of the pearling stability line, leading to the pearled morphology seen in Figure\,\ref{f:Jaylan} (center/b and c). Further increase in $\eta_d$ leads
to a defect laden filament structure, see Figure\,\ref{f:Jaylan} (center/d) consistent with slight crossing of the filament pearling stability line. This sequence, bilayer stability followed by
bilayer pearling and then filament pearling for $\mu_1>0$ and increasing $\eta_d$ is consistent with a value of $\xi$ in the range $[-0.4,-0.2]$ depicted in Figure\,\ref{f:Bdiag}.
 The qualitative agreement between the pearled morphologies observed in the FCH free energy and the self-consistent density-functional model is striking. Both models present
 emergent pearling with smaller, round holes, compare Figure\,\ref{f:Jaylan} (left/bottom-left) and (center/b), while fully emerged pearling leads to larger, pentagonal shaped holes,
compare Figure\,\ref{f:Jaylan} (left/bottom-right) and (center/c).

\subsection{Comparison of analytical bifurcation diagrams to experiments}

We compare the bifurcation structure derived for codimension 1 and 2 composite morphologies within the FCH gradient flow to experimental
results of Dicher and Eisenberg reprinted in Figure\,\ref{f:DE} and of Jain and Bates reprinted in Figure\,\ref{f:Dicher}.

\subsubsection{Bifurcations of Dicher and Eisenberg}

\begin{figure}[h!]

   \begin{tabular}{cc}
        \includegraphics[width=3.1in]{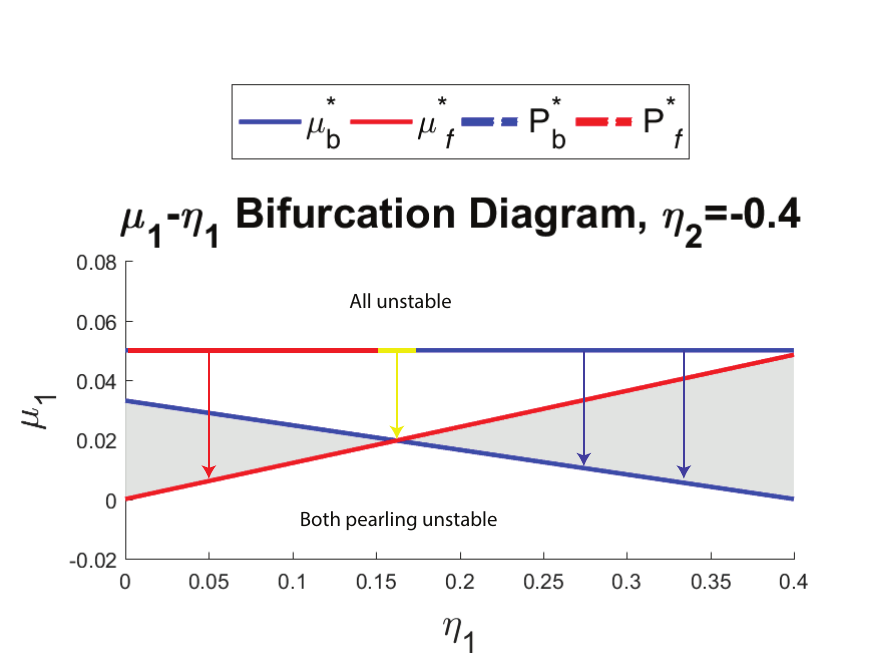} &
        \includegraphics[width=3.1in]{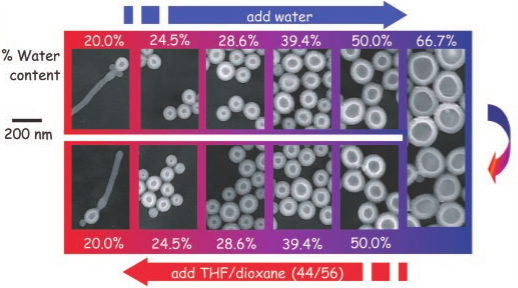}
       \end{tabular}

\caption{(left) A $\mu_1\!-\!\,\eta_1$ bifurcation diagram with $\xi=-0.5$ and $\eta_2=-0.4$ that shows the curve shortening lines for bilayers and filaments and
vertical arrows showing generic evolution of $\mu_1$ from initial data starting at $0.05$ for $\eta_1$ running from $0$ to $0.4$.
The color coding of the $\mu_1=0.05$ line indicates the final result of the end state with red denoting pure filament, yellow  coexistence of filament and bilayer, and blue
denoting pure bilayer; compare to Figure\,\ref{f:Dicher} (left) for increasing values of water in solvent phase. (right) Complete reversibility in a PS-PAA system under change in water solvent concentration, implying robust inhibition of pearling instabilities.  From \cite{discher2002polymer}. Reprinted with permission from AAAS.}
 \label{f:DE}
\end{figure}

The bifurcation experiments of Dicher and Eisenberg, \cite{discher2002polymer}, depicted in Figure\,\ref{f:Dicher} (left), characterize the end-states of casts
of PS-PAA amphiphilic diblocks dispersed in a water-dioxane solvent blend. Dioxane is a good solvent for both PEO and PS, but PS, like its commercial
relative styrofoam, is strongly hydrophobic in water. Increasing the water content from zero leads to end-state morphologies, in order:
solvability, only micelles, micelles and rods, only rods, rods and vesicles, and only vesicles. Here micelles are codimension 3, while rods and vesicles
refer to codimension 2 and 1 respectively. An increase in water content in the solvent is analogous to an increase in $\eta_1$ -- the energy release per unit of
interface formation, under constant diblock aspect ratio, hence constant $\eta_2.$  This morphological bifurcation sequence can be emulated within the FCH equation by fixing $
\xi=-0.5$ and $\eta_2=-0.4$ and allowing $\eta_1$ to vary from $0$ to $0.4$. The resulting $\mu_1\!-\!\,\eta_1$ bifurcation diagram is presented in Figure\,\ref{f:DE} (left). The
curve shortening lines are presented within the view, the pearling stability lines are outside the view.
For small values of $\eta_1$ filaments are dynamically favored, while bilayers are dynamically favored for larger values of $\eta_1$. In particular,
the horizontal line is color coded to depict a probable end-state morphology of initial data consisting of an admissible composite solution
with $\mu_1=0.05$. Here red denotes a pure filament end-state, blue a pure bilayer end-state and the yellow corresponds to a region of long-time
coexistence of the two morphologies due to the approximate equality of the critical values $\mu_b^*=\mu_\fil^*$ for $\eta_1=0.22.$ \emph{However} the bifurcation diagram lies within the pearling instability region of both bilayers
and filaments, and as a result the FCH predicts that neither of these pure states would persist. This is a limitation in the FCH model, with the well choice considered
here the intersection of the pearling bifurcation curves occurs at smaller values of $\eta_1$ than the intersection of the curve shortening lines.  Agreement with
these experimental results requires a robust inhibition of the pearling mechanism. Indeed, the experimental dynamics for the PS-PAA polymers considered here are largely reversible, as shown in Figure\,\ref{f:DE} (right), increasing and then decreasing the water content leads to a fully reversible sequence of morphological bifurcations. The presence of the pearling bifurcation generates complex morphologies with strong hysteresis, see Figure\,\ref{f:JB} (center). Matching this class of morphological bifurcation diagrams requires a tuning mechanism
within the well shape $W$ that affords robust pearling inhibition. Such a mechanism is proposed within the context of multicomponent models in \cite{promislow2017existence}
and \cite{DPV-19}.

\subsubsection{Bifurcations of Jain and Bates}

\begin{figure}[h!]
 \begin{tabular}{cp{3.5in}}
             \includegraphics[width=2.9in]{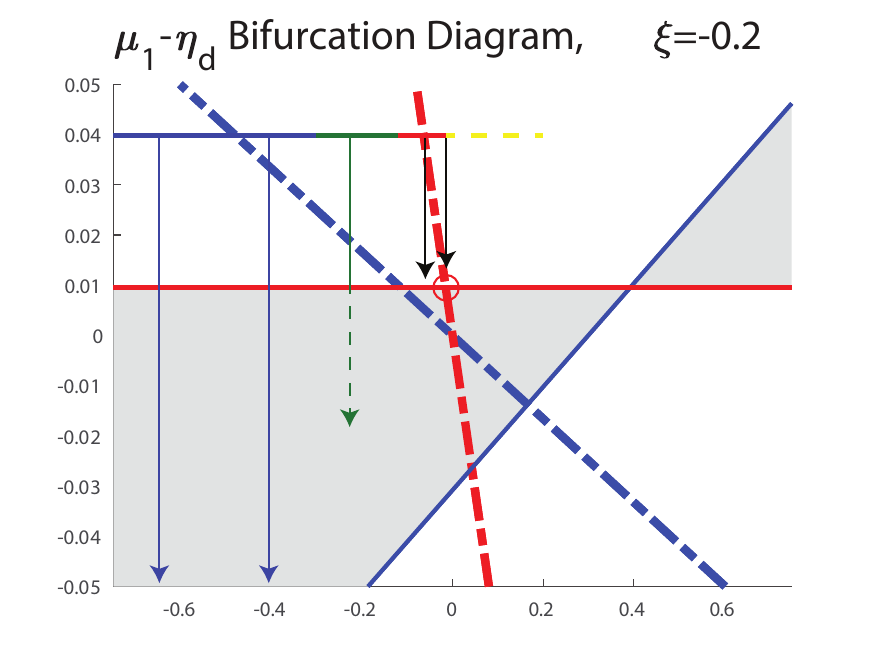}&\hspace{3.5in}
  \end{tabular}
  \vskip -2.0in
  \begin{tabular}{p{3.0in}cc}
          &     \includegraphics[width=1.1in]{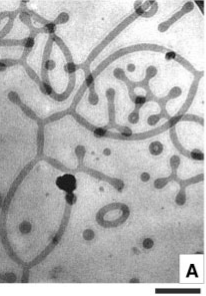}     &
               \includegraphics[width=2.05in]{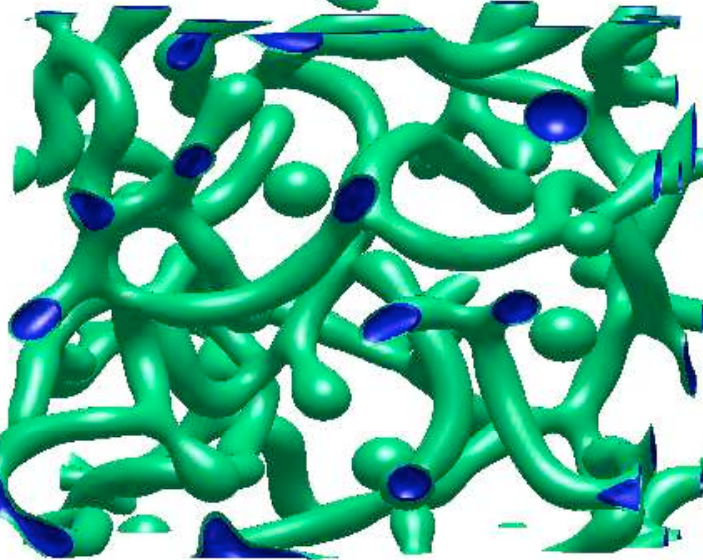}
   \end{tabular}
\caption{(left) A $\mu_1\!-\!\,\eta_d$ bifurcation diagram for $\xi=-0.2$ and
fixed $\eta_1=0.15$, showing the end states of initial data with $\mu_1=0.04$ and varying values of $\eta_2$, with color coding corresponding to probable end state: blue--bilayer,
green--network and defect structure (morphological complexit), red--filament, dotted yellow--filament and micelle.  Compare to experimental results of
Figure\,\ref{f:Dicher}(right) for $N_{PB}=170$, with increasing values of $w_{PEO}$ corresponding to  decreasing values of $\eta_2$. (center) Experimentally observed network, end cap, and $Y$-junction morphologies corresponding to the $C_Y$ phase of the bifurcation diagram from Figure\,\ref{f:Dicher} (right). Scale bar is $200$nm, From \cite{jain2003origins}. Reprinted with permission from AAAS. (right) End state of simulation of FCH gradient flow corresponding to green arrow, coarsened from random initial data, courtesy of Zhengfu Xu.}
 \label{f:JB}
\end{figure}

The experimental bifurcation diagram of Jain and Bates \cite{jain2003origins}, depicted in Figure\,\ref{f:Dicher} (right), shows the end-state morphology of dispersions of
PEO-PB with different polymer lengths and different weight fractions of PEO, $w_{\rm PEO}$, and hence different aspect ratios of the overall amphiphilic diblock.
For comparison Figure\,\ref{f:JB} (left) depicts the end states of
the FCH gradient flow  corresponding to values of $\xi=-0.2$ and $\eta_1=0.15$ for initial values of $\mu_1=0.075$ and various values of
$\eta_d$ arising from variation in $\eta_2$ which models the changes in diblock aspect ratio. Small, negative values of $\eta_2$ correspond to low values of
$w_{\rm PEO}$ while larger, positives values of $\eta_2$ correspond to larger values of $w_{\rm PEO}$ and to negative values
of $\eta_d$. The horizontal line in  Figure\,\ref{f:JB} (left) is color coded to depict a probable end-state of the FCH evolution starting from an admissible composite solution
with this value of $\mu_1$. On the left where $\eta_d<-0.5$,  an initial value of $\mu_1=0.075$, both codimension one and two morphologies are in their
 curve lengthening region and both would increase in surface area/length. However as $\mu_1$ is depleted, the suspension enters the curve shortening region for filaments,
 which will either vanish in finite time or pearl as $\mu_1$ crosses the red dotted line -- the end result is a pure bilayer state -- indicated by the blue color of the horizontal
line for this value of $\eta_d$.

We focus on the values of $\eta_d$ in $[-0.35, -0.05]$, for which the horizontal line is colored green, to indicate a region of \emph{morphological complexity}.
The initial data lies in the bilayer pearling region, but passes transiently through it to bilayer pearling stability region. Depending upon the form of the
initial bilayer morphology they may either fully pearl and form filament networks, or persist as bilayers and then grow after the return to pearling
stability. Filament networks formed from the pearling of a bilayer typically host many $Y$-junctions and end caps.  The filament network will expand until
 $\mu_1$ crosses the red-solid filament curve shortening line, at this point any defect free filaments will shrink, although the presence of any end-caps
and $Y$-junctions in a particular component will render its evolution unclear. This uncertainty is reflected in the dotted nature of the bottom half of the green vertical line.
The end result of the evolution is strongly dependent upon the form of the initial data, and will be very hysteretic in this regime. The uncertainty in the evolution is consistent
with co-existing bilayers, filament networks, $Y$-junctured filaments, and end-cap defects, see Figure\,{\ref{f:JB} (center) for a depiction of the experimental morphology
found in this regime and Figure\,\ref{f:JB} (right) for a corresponding end-state of simulation of the FCH gradient flow from random initial data.

For large values of $\eta_d$ the $\mu_1$ flow remains in the bilayer pearling instability region, the bilayers either will not form or will pearl and transform to
filament networks with  the end-result being a $Y$-juncture dominated filament network. The last prediction, for $\eta_d=-0.02$ corresponds to the triple
intersection (red circle) of the descending $\mu_1$ line, the filament pearling (red-dotted), and curve shortening (red-solid) lines.
The arrival of $\mu_1$ to the filament curve shortening line from above is consistent with a stable filament phase, but the emergent pearling bifurcation signals a
transition to end-cap and micelle formation, and corresponds to the possible coexistence of filament and micelle phases. This transition from filament to
filament-micelle is reflected in the dotted nature of the yellow $\mu_1$ horizontal line for $\eta_d>-0.02$. A final transition to a pure micelle stage is
plausible but is outside the scope of our analysis.  The overall trend depicted in Figure\,\ref{f:JB} (left) suggests increasing $\eta_d$ at fixed $\eta_1$ results in bilayers, bilayers mixed with defective filaments (end caps and $Y$-junctions), filaments, and filaments coexisting with micelles. This bifurcation sequence is in excellence qualitative agreement with
the $N_{\rm PB}=170$ bifurcation sequence depicted in Figure\,\ref{f:Dicher} (right) as the PEO weight fraction $w_{\rm PEO}$ {\sl decreases}  from high values to low values,
corresponding to decreasing $\eta_2$ and hence increasing $\eta_d$ subject to constant values of $\eta_1$.

The bifurcation sequence depicted in Figure\,\ref{f:Dicher} (right) for $N_{\rm PB}=45$ corresponds to a much shorter,
stiffer diblock polymer. In this regime the network and defect-laden filament phase
are not observed, rather increasing $w_{\rm PEO}$ weight fraction leads to the codimensional bifurcation sequence which leads from bilayers,
to coexistence of bilayers and filaments, to filaments, to coexistence of filaments and micelles, and finally to micelles.  While the general trend of the codimensional
bifurcation sequence is supported by the competitive geometric motion and its bifurcations, as in Figure\,\ref{f:DE} (left), we reiterate that
within the context of the scalar version of the FCH free energy presented herein, the pearling bifurcation cannot be fully suppressed.

\section{Discussion}
\label{s:Discussion}

We have presented a multiscale analysis of the $H^{-1}$ gradient flow of the FCH free energy corresponding to initial data close to dressings
of admissible codimension one (bilayer) and codimension two (filament) morphologies. We derived their curvature driven flow which couples to
the evolving, spatially constant far-field chemical potential, $\mu_1.$ This flow is the basis of the morphological competition, which barring
other bifurcations leads to an end state corresponding to the dynamically favored morphology with the lower critical value, $\mu_b^*$ or $\mu_\fil^*$,
of the far field chemical potential. In particular we identify regimes in which the geometric flow leads to growth or evanescence of each phase
through curve shortening or regularized curve lengthening. Combining the curve shortening/regularized curve lengthening bifurcation with the
pearling bifurcation results allows a characterization of the evolution of defect-free connected components of bilayer and filament morphologies.
Our analysis predicts that codimension one and two structures do not generically coexist on the long, $t=O(\eps^{-1})$ time scale within the strong FCH gradient flow.
Experimental results show transitions between distinct codimensional phases with relatively large margins of coexistence, however as remarked in the experimental literature,
the transients associated to these experiments are long compared to experimental patience and transients required to achieve a single phase of codimensional morphology
may require months to years, \cite{jain2004consequences}.

We compared the analytical bifurcation diagram to results from numerical simulations of the FCH gradient flow, to simulations of
a self-consistent mean-field density functional model for amphiphilic polymers, and to three sets of experimental studies of amphiphilic polymers.
We find that the self-consistent mean-field density functional model proposed in \cite{fraaije2003model} predicts a bifurcation sequence of
radial bilayers, pearled bilayers, and filaments with end-cap defects in strong qualitative agreement with the FCH bifurcation structure, in
particular the structure of the pearled spherical bilayers computed by both models are in excellent agreement. The experimental bifurcation analysis of Jain and Bates, \cite{jain2003origins},
was conducted at two polymer lengths, long polymers with $N_{\rm PB}=170$ and
shorter ones with $N_{\rm PB}=45.$ We find strong qualitative agreement between the FCH bifurcation structure and the experimental results for the longer chains, with the FCH results
suggesting that the development of \emph{morphological complexity} could be produced in a region in which
bilayer pearling and filament curve shortening bifurcations lie in close proximity, see section 4.4.2. The passage of the far-field chemical potential through these
bifurcations engenders network morphologies with $Y$-junctions, end-caps, and stable pearled filaments. The analytical basis of the complexity suggests a
mechanism for hysteresis: the passage through these sequences of bifurcations will not be readily reversed by a non-adiabatic return of the state variables.

The FCH model with the parameter choices presented herein does not qualitatively reproduce all of the experimental results. Both the solvent quality
bifurcation experiments of Dicher and Eisenberg, \cite{discher2002polymer} and the short polymer $N_{\rm PB}=45$ experiments
of Jain and Bates do not show evidence of pearling bifurcations. Complex morphology is not exhibited, and Dicher and Eisenberg
show that for their experimental parameters the morphology is remarkably reversible: hysteresis is not observed. However, the Dicher and Eisenberg
experimental bifurcation structure is well described by the dynamical favoritism arising from the morphological competition.  To match the full range of
experimental results the FCH free energy needs a tunable mechanism to robustly inhibit the pearling bifurcation. These can be achieved in two ways.
The first is physically motivated: pearling is a modulation of bilayer width, and shorter polymers are stiffer and lead to bilayers with a less compressible width.
The compressibility of the bilayer is tunable through the slope of the well $W$ at the value $u$ at which the bilayer density is greatest, generically the second $W=0$ crossing.
Tuning this slope to be large represents a stiffer diblock and may serve to inhibit the pearling mechanism.  The second mechanism is mathematically motivated: the pearling
bifurcation arises from a balance between the positive eigenvalue $\lambda_{b,0}$ of the linearization $L_{b,0}$, see (\ref{e:L0-def}) about the bilayer profile, and the negative
eigenspace of the Laplace-Beltrami operator, $\Delta_s$. The balance cannot occur if the operator  $L_{b,0}$ is non-self adjoint and whose positive real part
spectrum have non-zero imaginary parts. This arrangement can be tuned and detuned within a multicomponent model, providing precisely the desired
mechanism for robust pearling inhibition. This mechanism was discussed in section 5 of \cite{promislow2017existence} and forms the basis of the study of the
singularly perturbed systems in \cite{DPV-19}.

There remains a considerable amount to address within the family of FCH energies. While the FCH is at some level a phenomenological model with
parameter values that are not directly tunable from first principles, the FCH parameter values can be fit to experiments, or subscale molecular simulations,
much like Flory Huggins parameters.  Nonetheless, fitting the full form of a complicated well $W$, especially for multicomponent models, may be challenging.
Even within the simple model presented here, the role of micelles on morphological competition has not been addressed, and their stability, including their growth into
dumbbells and end-caped filaments is a primary instabilty mode that leads adiabatically from a codimension three structure into a defective codimension two structures.
 As there are almost no analytical characterization of defect modes, the micelle to dumbbell instability is worthy of study on a purely mathematical basis.

\section*{Acknowledgement}
The third author acknowledges supported by the National Science Foundation through grants DMS-1409940 and DMS-1813203. All authors thank Zhengfu Xu for
providing the numerical simulation which generated Figure\,\ref{f:JB} (right).
\newpage
\bibliographystyle{apalike}
\bibliography{NKbib}

\end{document}